\documentclass{gtpart}     
\usepackage{amscd}
\usepackage{tikz}
\tikzset{
    partial ellipse/.style args={#1:#2:#3}{
        insert path={+ (#1:#3) arc (#1:#2:#3)}
    }
}
\usetikzlibrary{matrix, calc}
\title{Tethers and homology stability for surfaces}

\author{Allen Hatcher}
\givenname{Allen}
\surname{Hatcher}
\address{Department of Mathematics, Cornell University, Ithaca, NY 14853, USA}
\email{hatcher@math.cornell.edu}
\urladdr{https://www.math.cornell.edu/~hatcher/}

\author{Karen Vogtmann}
\givenname{Karen}
\surname{Vogtmann}
\address{Mathematics Institute, University of Warwick, Coventry CV4 7AL, United Kingdom}
\email{kvogtmann@gmail.com}
\urladdr{http://www2.warwick.ac.uk/fac/sci/maths/people/staff/karen\_vogtmann/}
\arxivreference{1508.04334}
\arxivpassword{i5rbp}

\newtheorem{thm}{Theorem}[section]  
\newtheorem*{thm*}{Theorem}  
\theoremstyle{definition}
         
\newtheorem{proposition}{Proposition}[section] 
\newtheorem{corollary}[proposition]{Corollary}   
\newtheorem{example}[proposition]{Example}
\newtheorem{lemma}[proposition]{Lemma}
\newtheorem{remark}[proposition]{Remark}
\newcommand\Sub{S}
\newcommand\bdry{\partial}
\newcommand\iso{\cong}
\newcommand\incl{\hookrightarrow}
\newcommand\blackdot[2]{\fill[black](#1,#2) circle (.05)}
\newcommand\reddot[2]{\fill[red](#1,#2) circle (.05)}

\newcommand{\ord}[1]{\langle{#1}\rangle} 
\newcommand{\x}[1]{\textbf{#1}} 
\newcommand{\bs}[1]{{#1}^\prime} 
\DeclareMathOperator{\lk}{lk}
\DeclareMathOperator{\st}{st}

\begin{document}
 
\begin{abstract}
Homological stability for sequences $G_n \to G_{n+1} \to \cdots$ of groups is often proved by studying the spectral sequence associated to the action of $G_n$ on a highly-connected simplicial complex whose stabilizers are related to $ G_k$ for $k<n$.  When  $G_n$ is the mapping class group of a manifold, suitable simplicial complexes can be made using isotopy classes of various geometric objects in the manifold. In this paper we focus on the case of surfaces and show that by using more refined geometric objects consisting of certain configurations of curves with arcs that tether these curves to the boundary, 
the stabilizers can be greatly simplified and consequently also the spectral sequence argument.  We give a careful exposition of this program and its basic tools, then illustrate the method using braid groups before treating mapping class groups of orientable surfaces in full detail.
\end{abstract}

\maketitle


\section*{Introduction}
Many classical groups occur in sequences $G_n$ with natural inclusions $G_{n}\to G_{n+1}$.  Examples include the symmetric groups $\Sigma_n$, linear groups such as $GL_n$, the braid groups $B_n$, mapping class groups $M_{n,1}$ of surfaces with one boundary component, and automorphism groups of free groups $Aut(F_n)$.  A sequence of groups is said to be {\it homologically stable\/} if the natural inclusions  induce isomorphisms on homology 
$
H_i(G_n)\to H_i(G_{n+1})
$
for $n$ sufficiently large with respect to $i$. All of the sequences of groups mentioned above are homologically stable.  
This terminology is also slightly abused when there is no natural inclusion, such as for the mapping class groups $M_n$ of closed surfaces  and outer automorphism groups of free groups $Out(F_n)$; in this case, we say the series is homologically stable if the $i$th homology is independent of $n$, for $n$ sufficiently large with respect to $i$. 
   
Homology stability is a very useful property.  It sometimes allows one to deduce properties of the limit group $G_\infty = \lim_\to G_n$ from properties of the groups $G_n$; the classical example of this is Quillen's proof that various $K$-groups are finitely generated.   It is also useful in the opposite direction: it is sometimes possible to compute invariants of the limit group $G_\infty$, which by stability are invariants of the groups $G_n$; an example of this is the computation by Madsen and Weiss of the stable homology of the mapping class group.  Finally, there is the obvious advantage that homology computations   which are unmanageable for $n$ large can sometimes be done in $G_n$ for $n$ small.  

In unpublished work from the 1970's Quillen introduced a general method of proving stability theorems which was used by many authors in subsequent years (the earliest examples include \cite{Wagoner}, \cite{Vogtmann}, \cite{Charney},  \cite{vanderKallen}, and \cite{Har1}). The idea is to find a highly-connected complex $X_n$ on which  $G_n$ acts, such that stabilizers of simplices are isomorphic to $G_m$ for $m<n$. One then examines a slight variant of the equivariant homology spectral sequence for this action; this has  
\[
E^1_{p,q}=
\begin{cases}
H_q(G_n, \Z)&\text{for $p= -1$}\cr
\bigoplus_{\sigma\in \Sigma_p} H_q(stab(\sigma);\Z_{\sigma}) & \text{for  $p\geq 0$}
\end{cases}
\]
where $\Sigma_p$ is a set of representatives of orbits of $p$-simplices.  The fact that $X_n$ is highly connected implies that this spectral sequence converges to 0 for $p+q$ small compared to $n$, and the fact that simplex stabilizers are smaller groups $G_m$ means that the map $H_i(G_{n-1})\to H_i(G_n)$ induced by inclusion occurs as a $d^1$ map in the spectral sequence.  If one assumes the quotient $X_n/G_n$ is highly connected and one or two small conditions of a more technical nature are satisfied, then an induction argument on $i$ can be used to prove that this $d^1$ map is an isomorphism for $n$ and $i$ in the approximate range $n>2i$.

This is the ideal situation, but in practice the original proofs of homology stability were often more complicated because the complexes $X_n$ chosen had simplex stabilizers that were not exactly the groups $G_m$ for $m<n$.  
For the groups $Aut(F_n)$ and $Out(F_n)$ a way to avoid the extra complications was developed in the papers \cite{HV2} and \cite{HVW} with further refinements and extensions in \cite{HW}.  The idea was to use variations of the original complexes studied in earlier papers \cite{H2,HV1} that included more data.  In~\cite{HV2} this extra data consisted of supplementary $2$-spheres in the ambient $3$-manifold that were called ``enveloping spheres", while in~\cite{HVW} and~\cite{HW} this extra data was reformulated   in terms of arcs joining $2$-spheres to basepoints in the boundary of the manifold. These arcs could be interpreted as ``tethering" the spheres to the boundary.  

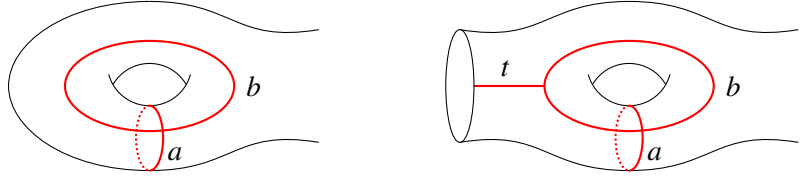
\begin{figure}
\begin{center}
  \begin{tikzpicture} [scale=.75]
   \draw (.5,0)  .. controls (.5,.75) and (1.5,1.5) ..  (3,1.5); 
    \draw (3,1.5)  .. controls (4.5,1.5) and (4.5,.75) ..  (6,1); 
    Ê \draw (.5,0)  .. controls (.5,-.75) and (1.5,-1.5) ..  (3,-1.5); 
     \draw (3,-1.5)  .. controls (4.5,-1.5) and (4.5,-.75) ..  (6,-1); 
\begin{scope}[xshift=3cm,yshift=.4cm]   
 \fill  [white] (-28:.75) arc (-28:-152:.75);
 \draw   (-15:.75) arc (-15:-165:.75);
 \draw [fill=white] (-30:.75) arc (30:150:.75);
 \end{scope}
 \draw [red, thick] (3,0) ellipse (1.5cm and .8cm);
 \draw[thick, red] (3,-.925) [partial ellipse=-90:90:.24cm and .575cm];
\draw[thick, red, densely dotted] (3,-.925) [partial ellipse=90:270:.24cm and .575cm];
\node (a) at (3.45,-1.2) {$a$};
\node (b)[right] at (4.5,0) {$b$};

\begin{scope}[xshift=8.5cm] 
   \draw (0,1)  .. controls (1.5,.75) and (1.5,1.5) ..  (3,1.5); 
    \draw (3,1.5)  .. controls (4.5,1.5) and (4.5,.75) ..  (6,1); 
    Ê \draw (0,-1)  .. controls (1.5,-.75) and (1.5,-1.5) ..  (3,-1.5); 
     \draw (3,-1.5)  .. controls (4.5,-1.5) and (4.5,-.75) ..  (6,-1); 
\begin{scope}[xshift=3cm,yshift=.4cm]   
 \fill  [white] (-28:.75) arc (-28:-152:.75);
 \draw   (-15:.75) arc (-15:-165:.75);
 \draw [fill=white] (-30:.75) arc (30:150:.75);
 \end{scope}
 \draw [red, thick] (3,0) ellipse (1.5cm and .8cm);
 \draw[thick, red] (3,-.925) [partial ellipse=-90:90:.24cm and .575cm];
\draw[thick, red, densely dotted] (3,-.925) [partial ellipse=90:270:.24cm and .575cm];
 \draw[thick, red] (.25,0) to (1.5,0);
  \draw  (0,0) ellipse (.25cm and 1cm);
\node (a) at (3.45,-1.2) {$a$};
\node (b)[right] at (4.5,0) {$b$};
\node (t) [above] at (.8,0) {$t$};
\end{scope}
 \end{tikzpicture}
 \caption{Chain on a closed surface and tethered chain on a surface with boundary.}
 \label{fig:chain}
 \end{center}
 \end{figure}

In the present paper we show how this tethering idea can be used in the case of mapping class groups of surfaces. 
As above, tethers are arcs to a point in the boundary, while at their other end they attach either to individual curves in the surface or to ordered pairs of curves intersecting transversely in one point.  We call such an ordered pair $(a,b)$ a {\it chain}    (see Figure~\ref{fig:chain}) and use the term ``curve" always to mean a simple closed curve.  

The classical curve complex $C(S)$ of a compact  orientable surface $S$ has vertices corresponding to isotopy classes of nontrivial curves in $S$ (where nontrivial means not bounding a disk or isotopic to a boundary component of $S$) and a set of vertices spans a simplex if the curves can be chosen to be disjoint.  Of particular interest is the subcomplex $C^0(S)$ formed by simplices corresponding to coconnected curve systems, that is, systems with connected complement.

In a similar way we can define a complex $Ch(S)$ of chains, with simplices corresponding to isotopy classes of systems of disjoint chains;  note that such a system is automatically coconnected, and there are no trivial chains to exclude.  When $S$ has nonempty boundary we can  
also form complexes $TC(S)$ and $TCh(S)$ of systems of disjoint tethered curves or tethered chains. In the case of $TC(S)$ we assume the curves without their tethers form coconnected curve systems, and in the case of $TCh(S)$ we assume the tether to a chain attaches to the $b$-curve.  A further variant that is useful for proving homology stability is a complex $DTC(S)$ of {\it double-tethered\/} curves, by which we mean curves with two tethers attached at the same point of the curve but on opposite sides, and curve systems are again assumed to be coconnected.  Note that shrinking the tether of a tethered chain to the point where it attaches to $\bdry S$ converts the $b$-curve to a double tether for the $a$-curve. More precise definitions for these complexes are given in Section~\ref{sec:tethers}, including extra data specifying where the tethers attach in $\bdry S$.

Our main new result is:

\begin{thm*}
If $S$ is a compact orientable surface of genus $g$ then the complexes $TC(S)$, $DTC(S)$, $Ch(S)$, and $TCh(S)$ are all $(g-3)/2$-connected.
\end{thm*}

Recall that a space $X$ is $r$-connected if $\pi_i(X)=0$ for $i\leq r$, which makes sense even if $r$ is not an integer.  Thus $r$-connected means the same as $\lfloor r\rfloor$-connected.  In particular $(-1)$-connected means nonempty (every map of $\bdry D^0 = \varnothing$ to $X$ extends to a map of $D^0$ to $X$) and $r$-connected for $r<-1$ is an empty condition.  

The mapping class group of $S$ acts on all these complexes.  A nice feature of the action on $TCh(S)$ is that the stabilizer of a vertex is exactly the mapping class group for a surface of genus one less but the same number of boundary components.  This is because cutting $S$ along a tethered chain reduces the genus by one without changing the number of boundary components.  Similarly, the stabilizer of a $k$-simplex is the mapping class group of a surface with genus reduced by $k+1$ and the same number of boundary components.  This makes  $TCh(S)$ ideal for the spectral sequence argument proving homology stability with respect to increasing genus with a fixed positive number of boundary components.  Actually it turns out to be slightly more efficient to use the complex $DTC(S)$, or a variant of it where the double tethers attach to basepoints in two different components of $\bdry S$ and the ordering of the tethers at these basepoints satisfies a compatibility condition. With this complex a single spectral sequence suffices to prove both that the homology stabilizes with respect to genus (namely, $H_i$ of the mapping class group is independent of $g$ for $g\geq 2i+2$) and that the stable homology does not depend on the number of boundary components as long as this number is positive.  In order to extend this to closed surfaces we need to work with a complex that does not involve tethers, and we use a version of $Ch(S)$ in which chains are oriented and systems of oriented chains are ordered. 
(Even if one is interested only in closed surfaces it is necessary to consider the case of non-empty boundary in order to have a way to compare mapping class groups in different genus.) 

The best stable dimension range that these simple sorts of spectral sequence arguments can yield has slope $2$, as in the inequality $g\geq 2i+2$.  This is not the optimal range, which has slope $3/2$, arising from more involved spectral sequence arguments.  See \cite{B,R-W,W1} for details.

The complexes of chains and tethered chains that we show are highly connected have found other recent applications as well in~\cite{PS} and~\cite{RWW}. 
In higher dimensions the natural analog of a tethered chain is a pair of $k$-spheres in a smooth manifold~$M^{2k}$ intersecting transversely in a single point, together with an arc tethering one of the spheres to a basepoint in~$\bdry M$. These tethered sphere-pairs play a central role in recent work of Galatius and Randal-Williams~\cite{GRW} on homology stability for $B{\rm Diff}(M)$ for certain $2k$-dimensional manifolds $M$ with $k>2$, including the base case that $M$ is obtained from a connected sum of copies of $S^k\times S^k$ by deleting the interior of a $2k$-ball.

Here is an outline of the paper.  In Section~\ref{sec:spectral}  we present the basic spectral sequence argument and in Section~\ref{sec:connectivity} we lay out the tools used to prove the key connectivity results.  In Section~\ref{sec:braids} we give a warm-up example illustrating the method in a particularly simple case, proving Arnold's homology stability theorem for braid groups.  In Section~\ref{sec:curvecomplexes} we give new, simpler proofs of results due to Harer about curve complexes and arc complexes  that will be used in Section~\ref{sec:tethers} to prove the main new connectivity statements.  Finally in Section~\ref{sec:MCGstability} we deduce homology stability for mapping class groups.

\medskip\noindent{\bf Remark.}  A draft version of this paper dating from 2006 and treating several other classes of groups has been informally circulated for a number of years. This current version focuses only on mapping class groups,  significantly simplifies several of the proofs in the earlier version, and also corrects a couple of errors. We thank Alexander Jasper for bringing one of these errors to our attention.

\section{The basic spectral sequence argument}\label{sec:spectral}

In this section we give the simplest form of the spectral sequence argument for proving homology stability of a sequence of group inclusions $\ \cdots \to G_n\to G_{n+1}\to G_{n+2}\to \cdots\ $. The input for the spectral sequence will be a simplicial action of $G_n$ on a simplicial or semi-simplicial complex $X_n$ for each $n$.  To deduce stability we will make the following assumptions which are stronger than necessary for stability but simplify the arguments and are satisfied in all but one of our applications.  The only exception arises in the proof of Theorem~\ref{thm:closedsurfaces}, where a short extra argument is required.
\begin{enumerate}
\item[(1)] $X_n$ has dimension $n-1$ and the action of $G_n$ is transitive on simplices of each dimension.
\item[(2)]  The stabilizer of a vertex is conjugate to $G_{n-1}$, and more generally the stabilizer of a $p$-simplex is conjugate to $G_{n-p-1}$. Moreover, the stabilizer of a simplex fixes the simplex pointwise. 
\item[(3)]  If $e$ is an edge of $X_n$ with vertices $v$ and $w$, then there is an element of $G_n$ taking $v$ to $w$ which commutes with all elements of the stabilizer of $e$.
\end{enumerate}
The dimension range in which homology stability holds will depend on the connectivity of $X_n,$ which must grow linearly with $n$. 
The best result that the method can yield is that $H_i(G_{n-1}) \to H_i(G_n)$ is an isomorphism for $n>2i+c$ and a surjection for $n=2i+c$, for some constant $c$.  
In the cases which occur in this paper we have the following  stable ranges.

\begin{thm}\label{thm:specseq}
Suppose the action of $G_n$ on $X_n$ satisfies conditions (1)--(3) for each $n$. Then:
\begin{enumerate}
\item[(a)]  If $X_n$ is $(n-3)$-connected for each $n$ then the stabilization $H_i(G_{n-1})\to H_i(G_n)$ is an isomorphism for $n>2i+1$ and a surjection for $n=2i+1$.
\item[(b)] If $X_n$ is $(n-3)/2$-connected for each $n$ then the stabilization $H_i(G_{n-1})\to H_i(G_n)$ is an isomorphism for $n>2i+2$ and a surjection for $n=2i+2$.
\end{enumerate}
\end{thm}

\begin{proof}
For $ G = G_n$ let $E_*G$ be a free resolution of $\Z$ by $\Z[G]$-modules, and let 
\[
\cdots \to C_p\to C_{p-1}\to\cdots\to C_0\to C_{-1}=\Z\to 0
\]
be the augmented simplicial chain complex of $X=X_n$.  The action of $G$ on $X$ makes $C_*$ into a complex of $\Z[G]$-modules, so we can take  the tensor product over $\Z[G]$ to form a double complex $C_*\otimes_G E_*G$.  Filtering this double complex horizontally and then vertically or vice versa gives rise to two spectral sequences, both converging to the same thing (see, e.g. \cite{Br}, VII.3).  

Using the horizontal filtration, the  $E^1_{p,q}$  term of the associated spectral sequence is formed by taking the $p$th homology of $C_*\otimes_G  E_qG$.  If we assume $X$ is highly connected, say $c(X)$-connected, then the complex $C_*$ is exact through dimension $c(X)$. Since $E_qG$ is free, ${C_*\otimes_G E_qG}$ is exact in the same range, so  $E^2_{p,q} = 0$ for $p\le c(X)$. In particular the spectral sequence converges to $0$ in the range $p+q \le c(X)$, so the same will be true for the other spectral sequence as well.

For the second spectral sequence, if we begin by filtering vertically instead of horizontally the associated $E^1_{p,q}$ term becomes $H_q(G; C_p)$.  By Shapiro's lemma (see, e.g., \cite{Br} p.~73) this reduces to 
\[
E^1_{p,q}=\bigoplus_{\sigma\in \Sigma_p} H_q(stab(\sigma);\Z_{\sigma})
\]
where $\Sigma_p$ is a set of representatives for orbits of $p$-simplices (if we consider a ``$(-1)$-simplex" to be empty, with stabilizer all of $G$), and  $\Z_\sigma$ is $\Z$ twisted by the orientation action of $stab(\sigma)$ on $\sigma$.  In our case we assume the action is transitive on $p$-simplices so there is only one term in the direct sum. We also assume $stab(\sigma)$ is conjugate to $G_{n-p-1}$ and fixes $\sigma$ pointwise, so $\Z_\sigma$ is an untwisted $\Z$ and the $E^1$ terms become simply  
$E^1_{p,q}=H_q(G_{n-p-1})$ with untwisted $\Z$ coefficients understood. 

Returning to the general case, the $q$-th row of the $E^1$ page is the augmented chain complex of the quotient $X/G$ with coefficients in the system  $\{H_q(stab(\sigma))\}$.  The $d^1$-differentials in this chain complex can be described explicitly as follows.  For a simplex $\sigma \in \Sigma_p$, the restriction of $d^1$ to the summand $H_q(stab(\sigma))$ will be the alternating sum of partial boundary maps $d^1_i\colon H_q(stab(\sigma))\to H_q(stab(\tau))$ where $\tau \in \Sigma_{p-1}$ is the orbit representative of the $i$th face $\partial_i\sigma$ and $d^1_i$ is induced by the inclusion $stab(\sigma) \to stab(\partial_i\sigma)$ followed by the conjugation that takes this stabilizer to $stab(\tau)$.

Homology stability is proved by induction on the homology dimension $i$, starting with the trivial case $i=0$. The sort of result we seek is that the stabilization map $H_i(G_{n-1})\to H_i(G_n)$ is an isomorphism for $n>\varphi(i)$ and a surjection for $n=\varphi(i)$, for a linear function $\varphi$ of positive slope.  

The map $d=d^1\colon E^1_{0,i}\to E^1_{-1,i}$ in the second spectral sequence constructed above is the map on homology induced by the inclusion of a vertex stabilizer into the whole group; by assumption this is the map $H_i(G_{n-1})\to H_i(G_n)$ induced by the standard inclusion $G_{n-1}\to G_n$, so this is the map we are trying to prove is an isomorphism. In the situation we are considering the $E^1$ page of the spectral sequence 
has the following form:
 \[
\begin{tikzpicture}
  \matrix (m) [matrix of math nodes,
                      nodes in empty cells, row sep=1em, column sep=4ex, minimum width=10ex]
  {
 H_i(G_n) &   H_i(G_{n-1}) &  H_i(G_{n-2}) & \cdots\hspace{3.5em}\\ 
\hspace{1em} \cdots\hspace{.5em} &  H_{i-1}(G_{n-1}) &    H_{i-1}(G_{n-2})   & H_{i-1}(G_{n-3})& \cdots\hspace{3.5em} \\ 
&  &   \hspace{3.5em}{\cdots} \hspace{.5em}   & H_{i-2}(G_{n-3}) &\cdots \hspace{3.5em}  \\ [3em]
 H_0(G_n) &   H_0(G_{n-1}) &  H_0(G_{n-2}) & H_0(G_{n-3})&\cdots\hspace{3.5em} \\ 
 p=\text{--}1 &  0 & 1 & 2&\cdots\hspace{3.5em} \\ 
};
 \path[-stealth]
    (m-1-2) edge node [above] {$d$} (m-1-1)
   (m-1-3) edge (m-1-2)
   (m-1-4) edge (m-1-3)
   (m-2-2) edge (m-2-1)
   (m-2-3) edge (m-2-2)
     (m-2-4) edge (m-2-3)
   (m-2-5) edge (m-2-4)
    (m-3-4) edge (m-3-3)
     (m-3-5) edge (m-3-4)
     (m-4-2) edge (m-4-1)
    (m-4-3) edge (m-4-2)
     (m-4-4) edge (m-4-3)
      (m-4-5) edge (m-4-4);
\node [left] (i) at (-6.5,2.55) {$i$};
\node [left] (ione) at (-6.5,1.55) {$i{-}1$};
\node [left] (itwo) at (-6.5,.55) {$i{-}2$};
\node [left] (q) at (-6.5,-1.55) {$q = 0$};
\draw (-6.45,-3) -- (-6.45,3);
\draw (-7.5,-2.2) -- (5.25,-2.2);
\end{tikzpicture}
\]
We first consider the argument for showing that $d$ is surjective. 
If $i$  is less than the connectivity $c(X)$ of $X$,
then the terms $E^\infty_{p,q}$ must be zero for $p+q\le i-1$, and in particular $E^{\infty}_{-1,i}$ must be zero. 
We will show that every differential $d^r$ with target $E^r_{-1,i}$  for $r>1$ is the zero map because its domain is the zero group, so the only differential that can do the job of killing $E^*_{-1,i}$ is $d$, which must therefore be onto. Thus it will suffice to show that $E^2_{p,q}=0$ for $p+q\le i$ and $q<i$. These groups are the reduced homology groups of $X/G$ with coefficients in the system of groups $\{H_q(stab(\sigma))\}$.  We will argue that these coefficient groups can be replaced by $H_q(G)$, with a suitable induction hypothesis, so that $E^2_{p,q}=\widetilde{H}_p(X/G;H_q(G))$, still assuming $p+q\le i$ and $q<i$. Thus the groups $E^2_{p,q}$ with $p+q\le i$ and $q<i$ will be zero once we know that the connectivity $c(X/G)$ is large enough, namely $c(X/G)\geq p$. Since we have $p+q\leq i$ and $q\geq 0$, the condition $c(X/G)\geq p$ can be reformulated as $c(X/G)\geq i$. 

As explained earlier, the $d^1$ differentials are built from maps induced by inclusion followed by conjugation. These maps fit into commutative diagrams
\[
\begin{tikzpicture}
\matrix (m) [matrix of math nodes, row sep=1.5em, column sep=2.5em, minimum width=2em]
{H_{q}(stab(\sigma)) & H_{q}(stab(\tau)) \\
H_{q}(G)& H_{q}(G) \\};
\path[-stealth]
   (m-1-1) edge (m-2-1) edge    (m-1-2)
   (m-2-1) edge  (m-2-2)
   (m-1-2) edge   (m-2-2);
\end{tikzpicture}
\]
where the vertical maps are induced by inclusion and the lower map is induced by conjugation in $G$, hence is the identity. If the vertical maps are isomorphisms we can then replace the coefficient groups in the $q$th row of the $E^1$ page by the constant groups $H_q(G)$. In our case with $E^1$ page displayed above, we would like the group $H_{i-1}(G_{n-3})$ and the groups to the left of it to be in the stable range, isomorphic to $H_{i-1}(G)$. Actually we can get by with slightly less, just having $H_{i-1}(G_{n-2})$ and the terms to the left of it isomorphic to the stable group and having $H_{i-1}(G_{n-3})$ mapping onto the stable group, since this is enough to guarantee that the homology of the chain complex at $H_{i-1}(G_{n-2})$ will be zero. Thus we want the relation $\varphi(i)\ge \varphi(i-1)+2$. The corresponding relation for smaller values of $i$ will take care of lower rows, by the same argument.

To summarize, we have shown that the stabilization $H_i(G_{n-1})\to H_i(G_n)$ will be surjective if $\varphi(i)\ge\varphi(i-1)+2$, assuming $i-1\le c(X_n)$ and $i\le c(X_n/G_n)$.

To prove that $d$ is injective the argument is similar, but with one extra step. If $i \le c(X_n)$ the term $E^\infty_{0,i}$ will be zero, and then it will suffice to show that all differentials with target  $E^r_{0,i}$ are zero, so the only way for $E^1_{0,i}$ to die is if $d$ is injective.  We can argue that the terms $E^2_{p,q}$ are zero for $p+q\le i+1$ and $q<i$ just as before, assuming again that $\varphi(i)\ge\varphi(i-1)+2$ but with the inequality $i\le c(X_n/G_n)$ replaced by $i+1\le c(X_n/G_n)$ since we have shifted one unit to the right in the spectral sequence.  The extra step we need for injectivity of $d$ is showing that the differential $d^1\colon E^1_{1,i}\to E^1_{0,i}$ is zero. This will follow from the assumption that for each edge $e$ of $X_n$ there is an element $g$ of $G_n$ taking one of the endpoints $v$ of $e$ to the other endpoint $w$, such that $g$ commutes with $stab(e)$. This guarantees that $d^1$ vanishes on the summand of $E^1_{1,i}$ corresponding to $e$ by our earlier description of $d^1$.  
Namely, conjugation by $g$ fixes $stab(e)$ and sends $stab(v)$ to $stab(w)$.  If $v_0$ is the vertex chosen to represent the vertex orbit and if $h_vv= v_0$ and $h_ww= v_0$, then the identifications of $stab(v)$ and $stab(w)$ with $stab(v_0)$ differ by conjugation by $h_wgh_v^{-1}$ so we have the following commutative diagram:  
\[
\begin{tikzpicture}
\matrix (m) [matrix of math nodes, row sep=.25em, column sep=3em, minimum width=2em]
{ & stab(v) &stab(v_0) \\
stab(e)&&\\
 & stab(w) &stab(v_0) \\};
\path[-stealth]
   (m-2-1) edge node[above] {$i$}   (m-1-2) edge  node[above] {$i$} (m-3-2) 
   (m-1-2) edge node[right] {$c_g$}  (m-3-2) edge node[above] {$c_{h_v}$}(m-1-3)
   (m-3-2) edge  node[above] {$c_{h_w}$} (m-3-3)
   (m-1-3) edge node[right] {$c_{h_wgh_v^{-1}}$}  (m-3-3)
;
\end{tikzpicture}
\]
Since $h_wgh_v^{-1}$ is an element of $stab(v_0)$,  conjugation by it induces the identity on $H_*(stab(v_0))$ and the previous diagram induces a commutative diagram
\[
\begin{tikzpicture}
\begin{scope}[yshift=-2.75cm]
\matrix (m) [matrix of math nodes, row sep=.5em, column sep=2.5em, minimum width=2em]
{ & H_{\ast}(stab(v)) &  \\
H_{\ast}(stab(e)) & &H_{\ast}(stab(v_0)) \\
  & H_{\ast}(stab(w)) &  \\};
\path[-stealth]
   (m-2-1) edge (m-1-2) edge   (m-3-2)
   (m-1-2) edge (m-2-3)  edge (m-3-2)
   (m-3-2) edge (m-2-3) 
  ;
  \end{scope}
\end{tikzpicture}
\]
Thus we see that the stabilization $H_i(G_{n-1})\to H_i(G_n)$ will be injective whenever $\varphi(i)\ge\varphi(i-1)+2$, assuming 
$i\le c(X_n)$ and $i+1\le c(X_n/G_n)$. 

The connectivity of the quotient $X_n/G_n$ is not hard to compute. Since we assume the action of $G_n$ is transitive on simplices of each dimension and the stabilizer of a simplex fixes it pointwise, $X_n/G_n$ can be identified with the quotient of the standard  simplex $\Delta^{n-1}$ obtained by identifying all of its $k$-dimensional faces for each $k$, where the identification preserves the ordering of the vertices. Thus $X_n/G_n$ is a semi-simplicial complex (or $\Delta$-complex) with one $k$-simplex for each $k\le n-1$.  It is easy to see that $X_n/G_n$ is simply-connected. Its simplicial chain complex has a copy of $\Z$ in each dimension $k\le n-1$ with boundary maps that are alternately zero and isomorphisms.  Therefore the reduced homology groups of $X_n/G_n$ are trivial below dimension $n-1$, while $H_{n-1}(X_n/G_n)$ is trivial when $n$ is odd and $\Z$ when $n$ is even.  Thus $X_n/G_n$ is $(n-2)$-connected

The condition $\varphi(i)\ge\varphi(i-1)+2$ is satisfied if we choose $\varphi(i)=2i+c$ for any constant~$c$.  It remains to determine $c$. 

To get surjectivity from the spectral sequence argument we need  $c(X_n/G_n)\geq i$ and $c(X_n)\geq i-1$. For injectivity we need one more degree of connectivity for each.   Consider first the inequalities involving $c(X_n/G_n)$.   We know $c(X_n/G_n)=n-2$ so we need $n\geq i+2$ for surjectivity and $n\geq i+3$ for injectivity.    We want surjectivity for $n\geq \varphi(i)=2i+c,$ for all $i\geq 1,$  (and injectivity for $n\geq 2i+c+1$) so any $c\geq 1$ works.  

 There remain the conditions $c(X_n)\geq i-1$ for surjectivity and $c(X_n)\geq i$ for injectivity.   In case (a) we have $c(X_n)=n-3$ giving the same $n\geq i+2$ for surjectivity and $n\geq i+3$ for injectivity as before, so $\varphi(i)=2i+1$ still works.    
For case (b) we have $c(X_n)=(n-3)/2$, so we need $n\geq 2i+1$ for surjectivity and $n\geq 2i+3$ for injectivity; in particular, taking $\varphi(i)=2i+2$ works for both. 
\end{proof}

\section{Connectivity tools}\label{sec:connectivity}

All of the complexes we will consider are of a certain type which we shall call, somewhat informally, {\it geometric complexes}.  Such a complex $X$ is a simplicial complex whose vertices correspond to isotopy classes of some type of nontrivial geometric object (for example arcs or curves in a surface, or combinations thereof), where {\it trivial\/} has different meanings in different contexts.  A collection of vertices $v_0,\cdots,v_k$ spans a $k$-simplex if representatives for the vertices can be chosen which are pairwise disjoint, and perhaps also satisfy some auxiliary conditions.  The corresponding set of isotopy classes defining the simplex of $X$ is also called a {\it system}, and the set of systems forms a partially ordered set (poset) $\widehat X$ under inclusion, whose geometric realization is the barycentric subdivision $\bs{X}$ of $X$. 

In this section we lay out a few general tools we will use for proving that various geometric complexes are highly connected. 

\subsection{Link arguments: rerouting disks to avoid bad simplices} \label{subsec:linksarguments}

We would like to relate $n$-connectedness of a simplicial complex $X$ to $n$-connectedness of a subcomplex $Y$.  We do this by finding conditions under which the relative homotopy groups $\pi_i(X,Y)$ are zero in some range $i\leq n$, so that the desired connectivity statements can be deduced from the long exact sequence of homotopy groups for $(X,Y)$.  Thus we wish to deform a map $f\colon (D^i,\bdry D^i)\to (X,Y)$ to have image in $Y$, staying fixed  on $\bdry D^i$.  We may assume $f$ is simplicial with respect to some triangulation of $D^i$, and then  
the idea is to deform $f$ by performing a sequence of local alterations in the open star of one simplex at a time, until $f$ is finally pushed into $Y$.  This method of improving the map is called a {\em link argument}.  We write $\lk(\sigma)$ for the link of a simplex $\sigma$ and $\st(\sigma)$ for the star; if the ambient complex $X$ needs to be specified we write $\lk_X(\sigma)$ and $\st_X(\sigma)$. 

We first identify a set of simplices in $X - Y$ as {\em bad} simplices, satisfying the following two conditions:
 \smallskip
\begin{itemize} 
\item[(1)]  Any simplex with no bad faces is in $Y$, where by a ``face" of a simplex we mean a subsimplex spanned by any nonempty subset of its vertices, proper or not.
\item [(2)]  If two faces of a simplex are both bad, then their join is also bad.
\end{itemize}
\smallskip
We call simplices with no bad faces  {\em good} simplices.  Bad simplices may have good faces, or faces which are neither good nor bad.  If $\sigma$ is a bad simplex we say  a simplex $\tau$ in $\lk(\sigma)$ is {\em good for $\sigma$} if any bad face of $\tau*\sigma$ is contained in $\sigma$.  The simplices which are good for $\sigma$ form a subcomplex of $\lk(\sigma)$ which we denote by $G_\sigma$.  

\begin{proposition}
\label{prop:link}
Let $X$, $Y$, and $G_\sigma$ be as above.  Suppose that for some integer $n\geq 0$ the subcomplex $G_\sigma$ of $X$ is $(n-\dim(\sigma)-1)$-connected for all bad simplices $\sigma$. Then the pair $(X,Y)$ is $n$-connected, i.e., $\pi_i(X,Y)=0$ for all $i\leq n$.
\end{proposition}

\begin{figure}
\begin{center}
 \begin{tikzpicture}[scale=.8] 
\begin{scope}[scale=1.5]
 \draw [fill=black!25] (0,-1) -- (0,1) -- (-1,0) -- (0,-1) -- (1,0) -- (0,1);
 \fill [red] (-1,0) circle (.04);
\fill [red] (1,0) circle (.04);
\fill [black] (0,1) circle (.04);
\fill [black] (0,-1) circle (.04);
 \node (mu) at (.13,0) {$\mu$};
 \draw [->] (2,0) to (3,0);
 \node (f) at (2.52,.21) {$f$};
 \end{scope}
\begin{scope}[xshift=8cm, yshift=-.6cm, scale=.6]
\draw [fill=red!10] (0,0) ellipse (4cm and 1cm);
\node (Y) at (0,-.5) {$Y$};
\draw [red] (-3,0) to (3,0);
\reddot{-3}{0};\reddot{-1}{0};\reddot{1}{0};\reddot{3}{0};
\draw (-1,3) to (1,3);
\blackdot{-1}{3}; \blackdot{1}{3};
\node (sigma) at (0,3.5) {$\sigma=f(\mu)$};
\draw [fill=black!60, fill opacity=0.7](-1,3) --(3,0)--(1,3)--(-1,3); 
\draw (-1,3) --(-3,0)--(1,3)--(-1,3);
\draw (-1,3) --(-1,0)--(1,3)--(-1,3);
\draw (-1,3) --(3,0)--(1,3)--(-1,3);
\draw (-1,3) --(1,0)--(1,3)--(-1,3);
\draw[fill=black!20, fill opacity=0.85](-1,3) --(-3,0)--(1,3)--(-1,3); 
\end{scope}
\end{tikzpicture}
 
\vskip 1cm

\begin{tikzpicture}[scale=.8]
\begin{scope}[scale=1.5]
 \draw [fill=black!25]  (0,1) -- (-1,0) -- (0,-1) -- (1,0) -- (0,1);
\fill [black] (0,1) circle (.04);
\fill [black] (0,-1) circle (.04);
 \draw [red, thick] (-1,0) -- (-.3,0) -- (.3,0) -- (1,0);
 \fill [red] (-1,0) circle (.04);
\fill [red] (1,0) circle (.04);
\fill [red] (-.3,0) circle (.04);
\fill [red] (.3,0) circle (.04);
\draw (0,-1) -- (-.3,0) -- (0,1) -- (.3,0) -- (0,-1);
\draw [->] (1.7,0) to (3.3,0);
\node (f ) at (2.5,.22) {$f|_{\partial\mu}*g$};
\end{scope}
\begin{scope}[xshift=8cm, yshift=-.6cm, scale=.6]
\draw [fill=red!20] (0,0) ellipse (4cm and 1cm);
\node (Y) at (0,-.5) {$Y$};
\draw (-1,3) to (1,3);
\blackdot{-1}{3}; \blackdot{1}{3};
\draw [fill=black!30, fill opacity=0.5](-1,3) --(-3,0)--(3,0)--(-1,3);
\draw  (-1,3) --(-1,0);
\draw (-1,3) --(1,0)--(1,3)--(-1,3);
\draw (-1,3) --(3,0)--(1,3)--(-1,3);
\draw[fill=black!50, fill opacity=0.7](1,3) --(-3,0)--(3,0)--(1,3);
\draw   (-1,0)--(1,3) -- (1,0);
\draw [red] (-3,0) to (3,0);
\reddot{-3}{0};
\reddot{-1}{0};
\reddot{1}{0};
\reddot{3}{0};
\end{scope}
 \end{tikzpicture}
 \caption{Retriangulation of  $\st(\mu)$ and new definition of $f$.}\label{link}
\end{center}
\end{figure}
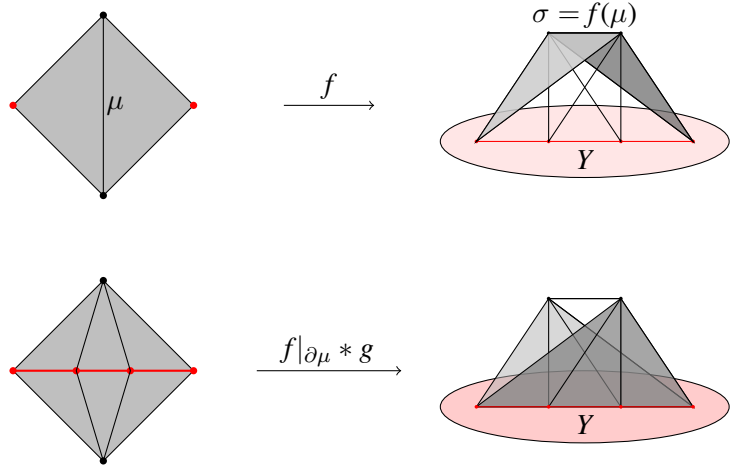
\begin{proof}  
We will show how to deform a map $f\colon (D^i,\bdry D^i)\to (X,Y)$  to have image in $Y$, staying fixed on $\bdry D^i$, provided that $i\leq n$. 
We may assume $f$ is simplicial with respect to some triangulation of $D^i$. Let $\mu$ be a maximal simplex of $D^i$ such that $\sigma= f(\mu)$ is bad (so in particular $\mu$ is not contained in $\bdry D^i$).  Then $f(\lk(\mu))\subset \lk(\sigma)$ is contained in $G_\sigma$, since otherwise there is some $\nu\in \lk(\mu)$ and face $\sigma_0$ of $\sigma$ such that $f(\nu)*\sigma_0$ is bad, in which case by property (2) $(f(\nu)*\sigma_0)*\sigma=f(\nu)*\sigma=f(\nu*\mu)$ is bad, contradicting maximality of $\mu$.   

We can assume the triangulation of $D^i$ gives the standard PL structure on $D^i$, so $\lk(\mu)$ is homeomorphic to $S^{i-k-1}$, where $k=\dim(\mu)\geq \dim(\sigma)$.  Since $G_\sigma$ is $(n-k-1)$-connected and $i\leq n$,  the restriction of $f$ to $\lk(\mu)$ can be extended to $g\colon D^{i-k}\to G_\sigma$, which we may take to be simplicial for some triangulation of $D^{i-k}$ extending the given triangulation on $S^{i-k-1}=\lk(\mu)$.  We retriangulate $\st(\mu)$ as $\bdry\mu*D^{i-k}$   
and redefine $f$ on this new triangulation to be $f |_{\bdry\mu}*g$ (see Figure~\ref{link}). 

The new map is homotopic to the old map, and agrees with the old map outside the interior of $\st(\mu)$, in particular on $\bdry D^i$.  Since simplices in $G_\sigma$ are good for $\sigma$, no simplices in the interior of $\st(\mu)$ have bad images. Since the original triangulation of $D^i$ was finite, in this way we can eventually eliminate all  $k$-simplices of $D^i$ with bad images without introducing simplices of higher dimension with bad images. Repeating the process in the new triangulation of $D^i$ for simplices of dimension $k-1$, then $k-2$, etc., we eventually eliminate all bad simplices from the image, so that by property (1) the image lies in $Y$.   
\end{proof}

 We give two applications of this proposition which we will use  in the rest of the paper.
\begin{corollary}
\label{lemma:subcomplex}
Let $Y$ be a subcomplex of a simplicial complex $X$, and suppose $X-Y$ has a set of bad simplices 
 satisfying (1) and (2) above. Then:
 \begin{itemize}
 \item [(a)] If $X$ is $n$-connected and $G_\sigma$ is $(n-\dim(\sigma))$-connected for all bad simplices $\sigma$, then $Y$ is $n$-connected. 
 \item[(b)] If \hskip2pt $Y$ is $n$-connected and $G_\sigma$ is $(n-\dim(\sigma)-1)$-connected for all bad simplices $\sigma$, then $X$ is $n$-connected. 
 \end{itemize}
\end{corollary}

\begin{proof}
Both statements follow from Proposition~\ref{prop:link} using the long exact sequence of homotopy groups for $(X,Y)$.  This is immediate for (b), while for (a) one should replace the $n$ in the proposition by $n+1$.
\end{proof}

Given any simplicial complex $X$ and a set of labels $S$, we can form a new simplicial complex $X^S$ whose simplices consist of the simplices of $X$ with vertices labeled by elements of $S$. Thus there are 
$|S|^{k+1}$ $k$-simplices of $X^S$ for each $k$-simplex of $X$.  

\begin{corollary}
\label{coloring}
Let $X$ be a simplicial complex and $S$ a set of labels.  
If $X$ is $n$-connected and the link of each $k$-simplex in $X$ is $(n-k-1)$-connected, then  $X^S$ is $n$-connected.  In the other direction, if $X^S$ is $n$-connected then so is $X$,   
without any condition on the links. 
\end{corollary}
\begin{proof}  
Choosing a label $s_0\in S$, we can regard $X$ as the subcomplex of $X^S$ consisting of simplices with all labels equal to $s_0$. There is then a retraction $r\colon X^S\to X$ which changes all labels to $s_0$. This implies the second statement of the corollary. For the first statement we will apply Corollary~\ref{lemma:subcomplex}(b).   
Call a simplex of $X^S$ {\em bad\/} if none of its vertex labels is equal to $s_0$.  It is immediate that the set of bad simplices satisfies (1) and (2).  If $\sigma$ is a bad simplex, then a simplex in $\lk(\sigma)$ is good for $\sigma$ if and only if all of its labels are $s_0$, so that $G_\sigma$ is isomorphic to $\lk_X(r(\sigma))$  and Corollary~\ref{lemma:subcomplex}(b)  
applies.      
\end{proof}

\begin{example}
{\rm If $X$ is the $p$-simplex $\Delta^p$, one might think the lemma could be applied for all $n$ to conclude that $(\Delta^p)^S$ was contractible. However, it can only be applied for $n\le p-1$ since for $n=p$ the hypothesis would say that the link of the whole simplex is $(-1)$-connected, i.e., nonempty, which is not the case. In fact $(\Delta^p)^S$ is the join of $p+1$ copies of the discrete set $S$, so it is $p$-dimensional and exactly $(p-1)$-connected if $S$ has more than one element.}
\end{example}

\subsection{Homotopy equivalence of posets}

The {\em geometric realization}   of a poset $P$ is the simplicial complex with one $k$-simplex for each totally ordered chain $p_0 < \cdots < p_k$ of $k+1$ elements $p_i\in P$.   An order-preserving map (poset map) between posets induces a simplicial map of their geometric realizations.  When we attribute some topological property to a poset or poset map we mean that the corresponding space or simplicial map has that property.  

For a poset map $\phi\colon P\to Q$ the fiber $\phi_{\leq q}$ over an element $q\in Q$ is defined to be the subposet of $P$ consisting of all $p\in P$ with $\phi(p)\leq q$. The fiber $\phi_{\geq q}$ is defined analogously.  The following  statement is known as Quillen's Fiber Lemma and is a special case of his Theorem A \cite{Q}.  We supply an elementary proof.   

\begin{proposition} \label{lemma:poset}
A poset map $\phi\colon P\to Q$ is a homotopy equivalence if all fibers $\phi_{\leq q}$ are contractible, or if all fibers $\phi_{\geq q}$ are contractible.
\end{proposition}

\begin{proof}  
There is no difference between the two cases, so let us assume the fibers $\phi_{\geq q}$ are contractible.  We construct a map $\psi\colon Q\to P$ inductively over the skeleta of $Q$ as follows. For a vertex $q_0$ we let $\psi(q_0)$ be any vertex in $\phi_{\geq q_0}$, which is non-empty since it is contractible. For an edge $q_0 < q_1$ both $\psi(q_0)$ and $\psi(q_1)$ then lie in $\phi_{\geq q_0}$ and we let $\psi$ map this edge to any path in $\phi_{\geq q_0}$ from $\psi(q_0)$ to $\psi(q_1)$. Extending $\psi$ over higher simplices $q_0< \cdots <q_k$ is done similarly, mapping them to $\phi_{\geq q_0}$ extending the previously constructed map on the boundary of the simplex.

We claim that $\psi$ is a homotopy inverse to $\phi$.  The composition $\phi\psi$ sends each simplex $q_0<\cdots < q_k$ to the subcomplex $Q_{\geq q_0}$. These subcomplexes are contractible, having minimum elements, so one can construct a homotopy from $\phi\psi$ to the identity inductively over skeleta of $Q$.  Similarly $\psi\phi$ is homotopic to the identity since it sends each simplex $p_0<\cdots < p_k$ to the contractible subcomplex $\phi_{\geq \phi(p_0)}$.
\end{proof}

We will often apply this proposition to the poset $\widehat X$ of simplices in some simplicial complex $X$. The geometric realization of this poset is the barycentric subdivision $\bs{X}$ of the complex. 
The following lemma characterizes the poset fibers.
 
\begin{lemma}\label{lemma:simplicialmap}
Let $f\colon X\to Y$ be a simplicial map of simplicial complexes,  $\widehat X$  the poset of simplices in $X$, $\widehat Y$ the poset of simplices in $Y$, and $\widehat f\colon \widehat X\to \widehat Y$ the induced poset map.  Then for each simplex $\sigma$ of $Y$ we have the following relationships:
\begin{enumerate}
\item[(i)] $\widehat f_{\leq \sigma}$ is homeomorphic to $f^{-1}(\sigma)$. 
\item[(ii)]$\widehat f_{\geq \sigma}$ is homotopy equivalent to $\widehat f^{-1}(\sigma)$.
\item[(iii)] $\widehat f^{-1}(\sigma)$ is homeomorphic to $f^{-1}(y)$ where $y$ is the barycenter of $\sigma$.
\end{enumerate}
\end{lemma}

\begin{proof}
Statement (i) is immediate from the definitions:  $\widehat f_{\leq \sigma}$ is the set of all simplices $\tau$ such that  $f(\tau)$ is a face of $\sigma$.  

On the other hand,  $\widehat f_{\geq \sigma}$   is the set of all simplices $\tau$ such that $f(\tau)$ has $\sigma$ as a face.  Since $f$ is a simplicial map, some face of $\tau$ maps to $\sigma$; let $\tau_\sigma$ be the (unique) maximal such face.  The map $\tau\mapsto \tau_\sigma$ is a poset map $\widehat f_{\geq \sigma}\to \widehat f^{-1}(\sigma)$ whose upper fibers are contractible, having unique minimal elements. Thus $\widehat f_{\geq \sigma}$ is homotopy equivalent to $\widehat f^{-1}(\sigma)$, giving statement~(ii). Part (iii) is clear from the definitions. 
\end{proof}

The following is an immediate consequence of Proposition~\ref{lemma:poset} and Lemma~\ref{lemma:simplicialmap}: 

\begin{corollary} \label{corollary:contractiblefibers}
Let $f\colon X\to Y$ be a simplicial map of simplicial complexes. If $f^{-1}(\sigma)$ is contractible for all simplices $\sigma$ or if $f^{-1}(y)$ is contractible for all barycenters $y$, then $f$ is a homotopy equivalence. \end{corollary}

\noindent{\bf Remark.}  
For a simplicial map, contractibility of the fibers over barycenters implies contractibility of all fibers since the fibers over an open simplex are all homeomorphic. Other types of maps for which contractibility of fibers implies homotopy equivalence or at least weak homotopy equivalence include fibrations, quasifibrations, and microfibrations (see \cite{Weiss} for the last case).  The corollary implies that simplicial maps with contractible fibers are quasifibrations, but they need not be fibrations or microfibrations, as shown by the simple example of vertical projection of the letter $\mathsf{L}$ onto its base segment.

\subsection{Fiber connectivity}

\begin{lemma}
\label{lemma:fiber}
Let $f\colon X\to Y$ be a simplicial map of simplicial complexes.  Suppose that $Y$ is $n$-connected and the  fibers $f^{-1}(y)$ over the barycenters $y$ of all $k$-simplices in $Y$ are $(n-k)$-connected.  Then $X$ is $n$-connected.
\end{lemma}

\begin{proof}
Given a map $g\colon S^i\to X$ which we can assume is simplicial, we want to extend this to a map $G\colon D^{i+1}\to X$ if $i\leq n$.  In order to do this, we first consider the composition $h=fg\colon S^i\to Y$.  Since $Y$ is $n$-connected, we can extend $h$ to a simplicial map $H\colon D^{i+1}\to Y$.  We will use $H$ to construct~$G$, which we will do inductively on the skeleta of the barycentric subdivision $\bs{D}$ of $D^{i+1}$. 

We begin by replacing all complexes and maps by the associated posets of simplices and poset maps:
\[
\begin{tikzpicture}
  \matrix (m) [matrix of math nodes,row sep=2em,column sep=3em,minimum width=2em]
  {
     \widehat S^i & \widehat X \\
     \widehat D^{i+1}& \widehat Y \\};
  \path[-stealth]
    (m-1-1) edge (m-2-1)
            edge  node [above] {$\widehat g$} (m-1-2)
    (m-2-1.east|-m-2-2) edge node [above] {$\widehat H$}
             (m-2-2)
    (m-1-2) edge node [right] {$\widehat f$} (m-2-2)
            ;
\end{tikzpicture}
\]
Let $\tau$ be a vertex of $\bs{D}$, so  $\tau$ can be viewed as  a simplex of $D^{i+1}$ or as an element of $\widehat D^{i+1}$. Since $H$ is simplicial, $\sigma=H(\tau)$ has dimension at most   $i+1\leq n+1$ in $Y$. By the hypothesis and Lemma~\ref{lemma:simplicialmap},   $\widehat f_{\geq\sigma}$ is at least $(-1)$-connected, i.e., it is non-empty, so choose $x\in \widehat f_{\geq\sigma}$ and set $G(\tau)=x$.  We can assume this agrees with the given $g$ for $\tau\in \bdry \bs{D}$.

Now assume we have defined $G$ on the $(k-1)$-skeleton of $\bs{D}$, and let $\tau_0<\cdots<\tau_k$ be a $k$-simplex of $\bs{D}$.  Let $\sigma_i=H(\tau_i)$, and note that $\widehat f_{\geq \sigma_j}\subset \widehat f_{\geq \sigma_0}$ for all $j$.  By construction, then, $G$ maps the boundary of the simplex to $\widehat f_{\geq \sigma_0}$.  Since $H$ 
 is a simplicial map, it can only decrease the dimension of a simplex, so $\dim(\sigma_0)\leq i+1-k\leq  n+1-k$, and consequently $\widehat f_{\geq \sigma_0}$ is at least $(k-1)$-connected.  Therefore we can extend $G$ over the interior of the $k$-simplex $\tau_0<\cdots<\tau_k$, agreeing with the given $G$ on $\bdry \bs{D}$.  This gives the induction step in the construction of $G$. 
\end{proof}

\subsection{Flowing into a subcomplex} 

In this section we abstract the essential features of a surgery technique from \cite{H1} for showing that certain complexes of arcs on a surface are contractible, in order to more conveniently apply the method to several different situations later in the paper.

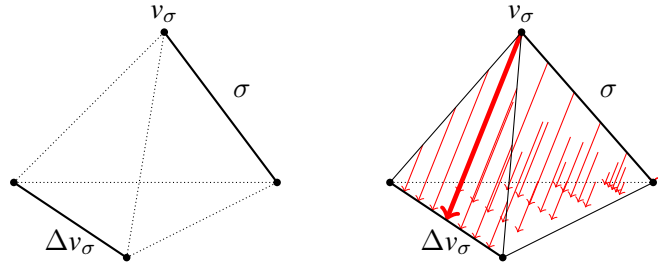
\begin{figure}
 \begin{center}
\begin{tikzpicture}[scale=.5]
 \coordinate (vs) at (4,4);
 \coordinate (w) at (7,0);
 \coordinate (x) at (0,0);
 \coordinate (y) at (3,-2);
  \draw [thick] (vs)-- (w);
  \fill [black] (w) circle (.1);
    \fill [black] (vs) circle (.1);
 \node [above] (vsig) at (vs) {$v_\sigma$};
  \node [above right] (sig) at ($(vs)!0.5!(w)$) {$\sigma$};
  \draw [thick]  (x) to (y);
   \draw [densely dotted]  (x) to (vs) to (y);
    \fill [black] (x) circle (.1);
   \fill [black] (y) circle (.1);
   \node [below] (Dv) at ($(x)!0.5!(y)$) {$\Delta v_\sigma$};
  \draw [densely dotted] (x) to (w);
  \draw  [densely dotted]  (y) to (w);
  \begin{scope}  [xshift=10cm]
   \coordinate (vs) at (3.5,4);
 \coordinate (w) at (7,0);
  \coordinate (x) at (0,0);
  \coordinate (y) at (3,-2);
  \draw [densely dotted] (x) to (w);
 \coordinate (m) at ($(x)!0.5!(y)$);
\draw [->,  red, line width=1.8pt] (vs)  to (m); 
\foreach \i in {1,...,4} {\draw [->, red] ($(vs)!0.2*\i!(w)$)  to ($(m)!0.2*\i!(w)$);} 
\foreach \i in {1,...,3} {\draw [->, red] ($(vs)!0.25*\i!(x)$)  to ($(m)!0.25*\i!(x)$);} 
\foreach \i in {1,...,3} {\draw [->, red] ($(vs)!0.25*\i!(y)$)  to ($(m)!0.25*\i!(y)$);} 
 \coordinate (r) at ($(m)!0.66!(y)$);
  \coordinate (s) at ($(vs)!0.66!(y)$);
  \foreach \i in {1,...,5} {\draw [->, red] ($(s)!0.2*\i!(w)$)  to ($(r)!0.2*\i!(w)$);} 
  \coordinate (r) at ($(m)!0.33!(y)$);
  \coordinate (s) at ($(vs)!0.33!(y)$);
   \foreach \i in {1,...,5} {\draw [->, red] ($(s)!0.2*\i!(w)$)  to ($(r)!0.2*\i!(w)$);} 
 \coordinate (r) at ($(m)!0.33!(x)$);
  \coordinate (s) at ($(vs)!0.33!(x)$);
   \foreach \i in {1,...,5} {\draw [->, red] ($(s)!0.2*\i!(w)$)  to ($(r)!0.2*\i!(w)$);} 
    \coordinate (r) at ($(m)!0.66!(x)$);
  \coordinate (s) at ($(vs)!0.66!(x)$);
   \foreach \i in {1,...,5} {\draw [->, red] ($(s)!0.2*\i!(w)$)  to ($(r)!0.2*\i!(w)$);} 
  \draw [thick] (vs)-- (w);
  \fill [black] (w) circle (.1);
\fill [black] (vs) circle (.1);
 \node [above] (vsig) at (vs) {$v_\sigma$};
  \node [above right] (sig) at ($(vs)!0.5!(w)$) {$\sigma$};
  \draw [thick]  (x) to (y);
   \draw   (x) to (vs) to (y);
    \fill [black] (x) circle (.1);
   \fill [black] (y) circle (.1);
   \node [below] (Dv) at ($(x)!0.5!(y)$) {$\Delta v_\sigma$};
  \draw  (y) to (w);
   \end{scope}
  \end{tikzpicture}
  \caption{A simplex $\sigma$, its preferred vertex $v_\sigma$, and flow lines in $\Delta v_\sigma * \sigma$}~\label{flow}
  \end{center}
  \end{figure}

Let $Y$ be a subcomplex of a simplicial complex $X$.
If $F\colon X\times I\to X$  is a deformation retraction into $Y$ then each $x\in X$ gives a path $F(x,t), 0\leq t\leq 1$  starting at $x$ and ending in $Y$.  In nice cases these paths fit together to give a flow on the complement of~$Y$.   What we want to do is to work backwards,  constructing a deformation retraction by first constructing a set of flow lines.    
Our flow lines will intersect each open simplex of $X$ which is not contained in $Y$ either transversely or in a family of parallel line segments.  
To specify these line segments, for each simplex $\sigma\in X-Y$  we choose   a preferred vertex $v=v_\sigma$  and a simplex $\Delta v$ in the link of $v$ in $X$ such that $\sigma * \Delta v$ is a simplex of $X$; then $\sigma*\Delta v$ is foliated by line segments parallel to the line from $v_\sigma$ to the barycenter of $\Delta v_\sigma$ (see Figure~\ref{flow}).
To show that the flow ends up in $Y$ we measure progress by means of a {\it complexity function\/} assigning a non-negative integer to each vertex of $X$ and taking strictly positive values on vertices not in $Y$.  This can be extended to be defined for all simplices of $X$, where the complexity of a simplex is the sum of the complexities of its vertices.

\begin{lemma}\label{lemma:flow} Let $Y$ be a subcomplex of a simplicial complex $X$ with a complexity function $c$ as above.  Suppose that for each vertex $v\in X - Y$  we have a rule for associating a simplex $\Delta v$ in the link of $v$ in $X$, and for each simplex $\sigma$ 
of $X$ not contained in $Y$ we have a rule for picking one of its vertices $v_\sigma \in X-Y$ so  that
\begin{itemize}
\item [{\rm (i)}]  the join $\sigma*\Delta v_\sigma$ is a simplex of $X$,
\item[{\rm (ii)}] $c(\Delta v)<c(v)$,
\item [{\rm (iii)}]  if $\tau$ is a face of $\sigma$ which contains $v_\sigma$, then $v_\tau=v_\sigma$.
\end{itemize}
Then $Y$ is a deformation retract of $X$.
\end{lemma} 

\begin{proof} 
For each simplex $\sigma$ not contained in $Y$ we construct flow lines in the simplex $\sigma*\Delta v_\sigma$ as described above, starting at $x\in \sigma$ and running parallel to the line from $v_\sigma$ to the barycenter of $\Delta v_\sigma$.  In terms of barycentric coordinates in the simplex $\sigma*\Delta v_\sigma$, viewed as weights on its vertices, we are shifting the weight on $v_\sigma$ to equally distributed weights on the vertices of $\Delta v_\sigma$, keeping the weights of other vertices fixed. When we follow the resulting flow on $\sigma*\Delta v_\sigma$, all points that actually move end up with smaller complexity by condition~(ii).  Thus after a finite number of such flows across simplices, each point of $\sigma$ follows a polygonal path ending in $Y$.  Condition~(iii) guarantees that the resulting flow is continuous on $X$, where we fix a standard Euclidean metric on each simplex and let each point flow at constant speed so as to reach $Y$ at time $1$. 
\end{proof}

\noindent {\bf Surgery flows.}  In this paper we will use Lemma~\ref{lemma:flow}  on various  complexes of arcs and curves on surfaces.  
The complexity function will count the number of nontrivial intersection points with a fixed arc, curve, or set of curves, 
and the simplex $\Delta v$ will be obtained using the surgery technique from \cite{H1} to decrease the number of intersection points.   The vertex $v_\sigma$ will be  an ``innermost" or ``outermost" arc or curve of $\sigma$, depending on the situation.  In order for this surgery process to be well-defined we must first put each arc or curve system $\{t_0,\cdots,t_k\}$  into normal form with respect to some fixed arc, curve, or curve system $t$, so that each $t_i$ has minimal possible intersection with $t$ in its isotopy class. In all cases we consider these normal forms are easily shown to exist; furthermore they are unique up to isotopy through normal forms, apart from the special situation that one $t_i$ is isotopic to $t$, in which case this $t_i$ can be isotoped across $t$ from one side to the other without always being in normal form during the isotopy.

\subsection{Ordered complexes}

In this subsection we prove a proposition that will be used at the very end of the paper when we extend the proof of homology stability for mapping class groups of non-closed surfaces to the case of closed surfaces.  We remark that this extension can also be proved without using the proposition, at the expense of complicating the spectral sequence argument and introducing an infinite-dimensional auxiliary complex.  This was the method used in three earlier papers in analogous situations:  \cite{H2} (page 53), \cite{HV2} (end of section 6), and \cite{HW} (proof of Theorem~5.1).

For a simplicial complex $X$ let $ \ord{X}$ be the ordered version of $X$, the semi-simplicial complex whose $k$-simplices are the $k$-simplices of $X$ with orderings of their vertices.  Thus there are $(k+1)!$ \hskip2pt $k$-simplices of $\ord{X}$ for each $k$-simplex of $X$.  For example if $X$ is a $1$-simplex then $\ord{X}$ has two vertices connected by two distinct edges.

Forgetting orderings gives a natural projection $\ord{X} \to X$.  This has cross-sections obtained by choosing an ordering of all the vertices of $X$ and using this to order the vertices of each simplex of $X$.  Thus $X$ is a retract of $\ord{X}$, so high connectivity of $\ord{X}$ implies the same high connectivity for $X$. We will be interested in the converse question of when high connectivity of $X$ implies high connectivity of $\ord{X}$. It is clear that $\ord{X}$ is connected if $X$ is, but the example of a $1$-simplex shows that this does not extend to $1$-connectedness. We therefore need some conditions on $X$, conditions that will be satisfied in our application. 
 
Generalizing a property of spheres with PL triangulations, a simplicial complex of dimension $n$ is called {\em Cohen-Macaulay } if it is $(n-1)$-connected and the link of each of its $k$-simplices is $(n-k-2)$-connected. If we drop the condition that $n$ is the dimension of the complex and only require it to be $(n-1)$-connected with the link of each $k$-simplex $(n-k-2)$-connected, then we have the notion of {\em weakly Cohen-Macaulay (wCM) of level~$n$}. In the existing literature, e.g., \cite{HW}, the term ``level"  is replaced by ``dimension", although this may be misleading since there is no restriction on the actual dimension.  

A simple observation is that a complex $X$ is wCM of level $n$ if and only its $n$-skeleton is wCM of level $n$. This is because the homotopy groups $\pi_i(X)$ for $i\leq n-1$ depend only on the $n$-skeleton, and a similar statement holds also for links in $X$. Thus $X$ being wCM of level $n$ is equivalent to its $n$-skeleton being Cohen-Macaulay (of dimension $n$).

Note that a complex which is wCM of level $n$ is automatically wCM of level $m$ for each $m<n$.  There is no need to require $n$ to be an integer, but if it is not, then wCM of level $n$ is the same as wCM of level $\lfloor n\rfloor$, the greatest integer $\leq n$, so allowing $n$ to be non-integral is just a matter of convenience.

Another observation is that if $X$ is wCM of level $n$ and $\sigma$ is a $k$-simplex in $X$ then the link of $\sigma$ is wCM of level $n-k-1$.  This is because for $\tau$ an $l$-simplex in $\lk(\sigma)$ we have $\lk_{\lk(\sigma)}(\tau)=\lk_X(\sigma*\tau)$ with $\sigma*\tau$ having dimension $k+l+1$, so $\lk_X(\sigma*\tau)$ has connectivity $n-(k+l+1)-2=(n-k-1)-l-2$.

\begin{proposition}\label{prop:OrderedComplex}  If a simplicial complex $X$ is wCM of level $n$ then the ordered complex $\ord{X}$ is $(n -1)$-connected.
\end{proposition}

In addition to the proof given below we give a different proof in the Appendix, following instead the approach in Proposition~2.14 of~\cite{RWW}. We include both proofs since they are of similar length and each has its own virtues.

 \begin{proof} 
As notation, we will use lowercase Greek letters for simplices of $X$, while simplices in $\ord{X}$ will be written as the ordered string of their vertices, e.g., $x_0x_1\cdots x_k$, sometimes abbreviated to $\x{x}=x_0x_1\cdots x_k$.

By induction on $n$ it suffices to show that $\pi_{n-1}(\ord{X})=0$. 
Given a map $f\colon \bdry D^n\to \ord{X}$, compose it with the projection $\ord{X}\to X$ to get a map $\bdry D^n\to X$.  Since $X$ is wCM of level $n$, this can be extended to a map $D^n\to X$.
Composing this extension with a section gives a map $ D^n\to\ord{X}$, whose restriction $g\colon \bdry D^n\to \ord{X}$ has the same projection as $f$.
Since $g$ is homotopically trivial, it suffices to show $f$ is homotopic to~$g$. 

In order to construct a homotopy we cover $X$ by the stars  $\st_X(\sigma)$ of its simplices and consider the corresponding cover of $\ord{X}$ by the ordered complexes $\ord{\st_X(\sigma)}$. We first show that each of these is $(n-1)$-connected, and then use this fact to build the homotopy between $f$ and $g$.  

\medskip\noindent {\em Claim.} {\em For each $k$-simplex $\sigma$ in $X$, $\ord{\st_X(\sigma)}$ is $(n-1)$-connected.}

\begin{proof} We may assume $X$ has dimension $n$ since higher-dimensional simplices have no effect on the relevant connectivities, as noted earlier.  Choose a vertex $a\in \sigma$. Then $\st_X(\sigma)$ is the cone $a*Y$ on $Y=\tau*\lk_X(\sigma)$, where $\tau$ is the face of $\sigma$ opposite $a$.  Since $\tau$ is wCM of level $k-1$ and $\lk_X(\sigma)$ is wCM of level $n-k-1$ as noted earlier, it follows that $Y$ is wCM of level $n-1$.  By induction the proposition is therefore true for $Y$ (and all of its links).  

Filter $\ord{\st_X(\sigma)}$ by subcomplexes $\Delta_i$, where $\Delta_i$ is the union of all ordered $n$-simplices $x_0\cdots x_{n}$ in $\ord{\st_X(\sigma)}$ with $a=x_j$ for some $j\leq i$. The lower-dimensional simplices of $\Delta_i$ thus have the form $x_0\cdots x_m$ with either   $x_j=a$ for some $j\leq i$ or no $x_j=a$. 
We will show that each $\Delta_i$ is $(n-1)$-connected by induction on $i$.   The first subcomplex $\Delta_0$ is the union of all ordered simplices of the form $ay_1\cdots y_n$ with $y_1\cdots y_n\in \ord{Y}$, i.e., it is the cone $a*\ord{Y}$ so is contractible.  Each $\Delta_i$ for $i>0$ is obtained from $\Delta_{i-1}$ by attaching all $n$-simplices of $\ord{\st_X(\sigma)}$ of the form $x_0\cdots x_{i-1}ay_{i+1}\cdots y_n$.  If we fix the ordered simplex  $\x{x}=x_0\cdots x_{i-1}$ and let the $y_j$  
vary we obtain a subcomplex  $\Delta_i(\x{x})$ of $\Delta _i$, i.e., $\Delta_i(\x{x})$ is the union of the ordered $n$-simplices of $\ord{\st_X(\sigma)}$ starting with $\x{x}a$.   

The subcomplex $\Delta_i(\x{x})$ decomposes as the join $\x{x}a*\ord{\lk_Y(\eta)}$, where $\eta$ is the (unordered) projection of $\x{x}$ to $X$.  In particular $\Delta_i(\x{x})$ is contractible since $\x{x}a$ is contractible.  
The intersection of $\Delta_i(\x{x})$ with $\Delta_{i-1}$ is $\bdry(\x{x}a) * \ord{\lk_Y(\eta)}$ since the only way a face of a simplex $x_0\cdots x_{i-1}ay_{i+1}\cdots y_n$ can lie in $\Delta _{i-1}$ is if at least one of the vertices $x_0,\cdots,x_{i-1},a$ is deleted.   Since $\bdry(\x{x}a)$ is an $(i-1)$-sphere, $\bdry(\x{x}a) * \ord{\lk_Y(\eta)}$ is the $i$-fold suspension of $\ord{\lk_Y(\eta)}$ (since join with $S^0$ is suspension and join is associative) so the connectivity of $\bdry(\x{x}a) * \ord{\lk_Y(\eta)}$ is $i$ more than the connectivity of $ \ord{\lk_Y(\eta)}$. By induction $\ord{\lk_Y(\eta)}$ is $((n-i-1)-1)$-connected, so $\bdry(\x{x}a) * \ord{\lk_Y(\eta)}$ is $(n-2)$-connected.
 An application of the Mayer-Vietoris sequence and the van Kampen theorem then shows that $\Delta_{i}(\x{x})\cup \Delta_{i-1}$ is $(n-1)$-connected.  

If we fix an ordered simplex $\x{z}=z_0\cdots z_{i-1}$   different from $\x{x}$, then the intersection of $\Delta_i(\x{z})$ and $\Delta_i(\x{x})$ is contained in $\Delta_{i-1}$ since simplices in the intersection can only be obtained by deleting at least one vertex of $\x{x}$ and of $\x{z}$ (and possibly other vertices).  We can then apply the above argument inductively to show that attaching finitely many complexes $\Delta_i(\x{x})$ to $\Delta_{i-1}$ preserves the connectivity $n-1$.  Since homotopy groups commute with direct limits, it follows that the entire subspace $\Delta_i$ is $(n-1)$-connected. Since $\Delta_n=\ord{\st_X(\sigma)}$ the claim is established.
\end{proof}

We  now proceed to build our homotopy $f\simeq g$.  The semi-simplicial complex $\ord{X}$ has the property that the vertices of each simplex are all distinct, so its barycentric subdivision $\bs{\ord{X}}$ is a simplicial complex.  We view $f$ and $g$ as maps $S^{n-1}\to\bs{\ord{X}}$ which we may take to be simplicial with respect to some triangulation of $S^{n-1}$.  We build the homotopy inductively on the skeleta of $S^{n-1}$.  

If $v$ is a vertex of $S^{n-1}$ then $f(v)$ and $g(v)$ project to the same vertex of $\bs{X}$, i.e., to the barycenter of some simplex $\sigma$ of $X$.  Since $\ord{\st_X(\sigma)}$ is $(n-1)$-connected 
there is a path in $\ord{\st_X(\sigma)}$  connecting $f(v)$ to $g(v)$, and we use this to define our homotopy on $v\times I$.  

Now let $s$ be any simplex of $S^{n-1}$
and assume we have already defined a homotopy $f\simeq g$ on $\bdry s$.    The projection of $f(s)$ (and hence $g(s)$) to $X$ is a simplex of $\bs{X}$, i.e., a chain $\sigma_0\subset\cdots\subset\sigma_k$ of simplices of $X$.   The stars of these simplices satisfy the reverse inclusions 
$\st_X(\sigma_0)\supset \cdots \supset \st_X(\sigma_k),$ hence the same is true for the ordered versions of these stars. We may assume by induction that the homotopy from $f$ to $g$ on $\bdry s$ takes place in $\ord{\st_X(\sigma_0)}$.  Since $\ord{\st_X(\sigma_0)}$ is $(n-1)$-connected and $\dim(\bdry s)\leq n-2$ we can extend the homotopy over the interior of $s$ so that its image lies in $\ord{\st_X(\sigma_0)}$.  This finishes the induction step.
\end{proof}

\section{A simple example: the braid group}\label{sec:braids}

As a warm up for our main case of mapping class groups let us first show how the method described in this paper 
gives a simple proof of homology stability for the classical braid groups $B_n$, where we are viewing $B_n$  as the mapping class group of an $n$-punctured disk. 

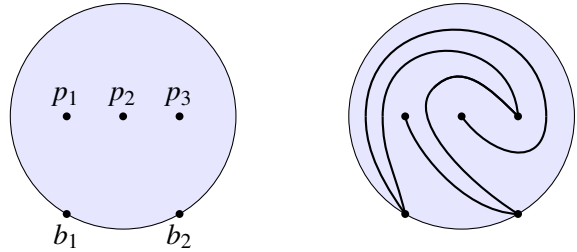
\begin{figure}
\begin{center}
\begin{tikzpicture} [scale=.75]

 \begin{scope}[xshift = -6cm] 
\draw [fill=blue!10] (0,0) circle (2);
\fill [black] (-1,0) circle (.07);
\node [above] (p2) at (-1,0) {$ {p_1}$};
\fill [black] (0,0) circle (.07);
\node [above] (p1) at (0,0) {$ {p_2}$};
\fill [black] (1,0) circle (.07);
\node [above] (p3) at (1,0) {$ {p_3}$};
 \fill [black] (-60:2) circle (.07);
 \node [below] (b1) at (-120:2) {$ {b_1}$};
  \node [below] (b2) at (-60:2) {$ {b_2}$};
\fill [black] (-120:2) circle (.07);
\end{scope}

\draw [fill=blue!10] (0,0) circle (2);
\fill [black] (-1,0) circle (.07);
\fill [black] (0,0) circle (.07);
\fill [black] (1,0) circle (.07);
 \fill [black] (-60:2) circle (.07);
\fill [black] (-120:2) circle (.07);

\draw [thick] (-60:2) .. controls (0,-1.1) and (-1,0) .. (-.5,.5);
\draw [thick] (1,0) .. controls (.5,.5) and (0,1) .. (-.5,.5);
\draw [thick] (-60:2) .. controls (0,-1.7) and (-.9,-.5) .. (-1,0);
\draw [thick] (1,0) .. controls (.5,.5) and (0,1) .. (-.5,.5);
\draw [thick] (-120:2) .. controls (-1.2,-.9) and (-1.4,-.7) .. (-1.4,0);
\draw [thick] (1,0) .. controls (.9,1.5) and (-1.4,1.6) .. (-1.4,0);
\draw [thick] (-120:2) .. controls (-1.6,-.9) and (-1.7,-.4) .. (-1.7,0);
\draw [thick] (-1.7,0) .. controls (-1.7,2.1) and (1.5,2) .. (1.5,0);
\draw [thick] (0,0) .. controls (.5,-.7) and (1.5,-1) .. (1.5,0);

 \end{tikzpicture}
 \caption{A system of tethers in a 3-punctured disk with two boundary points}\label{fig:tetherexample}
 \end{center}
  \end{figure}

We start by constructing a suitable ``tethered" complex with a $B_n$-action.  In fact, in this case the tethers will be all there is to the complex.  Consider a fixed disk $D$ with $d$ distinguished points $b_1,\cdots, b_d$ on the boundary and $n$ marked points or punctures $p_1,\cdots, p_n$ in the interior.  A {\it tether\/} is an arc in $D$ connecting some $p_i$ to some $b_j$ and disjoint from the other $p_k$'s and $b_k$'s. 
A {\it system of tethers\/} is a collection of tethers which are disjoint except at their endpoints, and with no two of the tethers isotopic.  See Figure~\ref{fig:tetherexample} for an example.

Define the tether complex $T=T_{n,d}$ to be the geometric
complex having one $k$-simplex for each isotopy class of systems of $k+1$ tethers, where the face relation is given by omitting tethers.  

\begin{proposition} \label{prop:surgery}
$T$ is contractible.
\end{proposition}

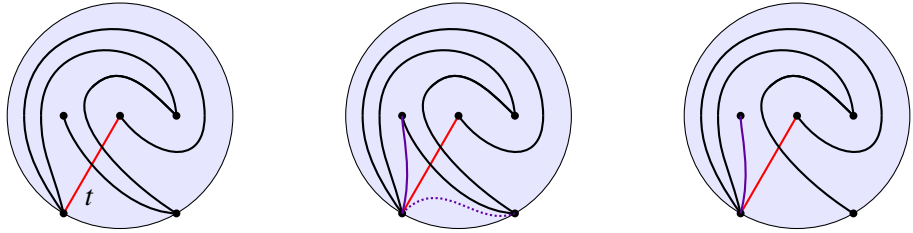
\begin{figure}
\begin{center}
 \begin{tikzpicture}  [scale=.75]
\draw [fill=blue!10] (0,0) circle (2);

\draw [red, thick] (-120:2) to (0,0);
\node (t) at (-.55,-1.45) {$t$};

\fill [black] (-1,0) circle (.07);
\fill [black] (0,0) circle (.07);
\fill [black] (1,0) circle (.07);
 \fill [black] (-60:2) circle (.07);
\fill [black] (-120:2) circle (.07);

\draw [thick] (-60:2) .. controls (0,-1.1) and (-1,0) .. (-.5,.5);
\draw [thick] (1,0) .. controls (.5,.5) and (0,1) .. (-.5,.5);
\draw [thick] (-60:2) .. controls (0,-1.7) and (-.9,-.5) .. (-1,0);
\draw [thick] (1,0) .. controls (.5,.5) and (0,1) .. (-.5,.5);
\draw [thick] (-120:2) .. controls (-1.2,-.9) and (-1.4,-.7) .. (-1.4,0);
\draw [thick] (1,0) .. controls (.9,1.5) and (-1.4,1.6) .. (-1.4,0);
\draw [thick] (-120:2) .. controls (-1.6,-.9) and (-1.7,-.4) .. (-1.7,0);
\draw [thick] (-1.7,0) .. controls (-1.7,2.1) and (1.5,2) .. (1.5,0);
\draw [thick] (0,0) .. controls (.5,-.7) and (1.5,-1) .. (1.5,0);


 \begin{scope}[xshift = 6cm] 

\draw [fill=blue!10] (0,0) circle (2);

\draw [red, thick] (-120:2) to (0,0);

\fill [black] (-1,0) circle (.07);
\fill [black] (0,0) circle (.07);
\fill [black] (1,0) circle (.07);
 \fill [black] (-60:2) circle (.07);
\fill [black] (-120:2) circle (.07);

\draw [thick] (-60:2) .. controls (0,-1.1) and (-1,0) .. (-.5,.5);
\draw [thick] (1,0) .. controls (.5,.5) and (0,1) .. (-.5,.5);
\draw [thick] (-60:2) .. controls (0,-1.7) and (-.9,-.5) .. (-1,0);
\draw [thick] (1,0) .. controls (.5,.5) and (0,1) .. (-.5,.5);
\draw [thick] (-120:2) .. controls (-1.2,-.9) and (-1.4,-.7) .. (-1.4,0);
\draw [thick] (1,0) .. controls (.9,1.5) and (-1.4,1.6) .. (-1.4,0);
\draw [thick] (-120:2) .. controls (-1.6,-.9) and (-1.7,-.4) .. (-1.7,0);
\draw [thick] (-1.7,0) .. controls (-1.7,2.1) and (1.5,2) .. (1.5,0);
\draw [thick] (0,0) .. controls (.5,-.7) and (1.5,-1) .. (1.5,0);

\draw [purple!50!blue, thick] (-120:2) .. controls (-.8,-1) and (-1,-.3) .. (-1,0);
\draw [purple!50!blue, thick, densely dotted] (-120:2) .. controls (-.3,-1) and (.3,-2) .. (-60:2);

 \end{scope}
  \begin{scope}[xshift = 12cm] 
\draw [fill=blue!10] (0,0) circle (2);

\draw [red, thick] (-120:2) to (0,0);

\fill [black] (-1,0) circle (.07);
\fill [black] (0,0) circle (.07);
\fill [black] (1,0) circle (.07);
 \fill [black] (-60:2) circle (.07);
\fill [black] (-120:2) circle (.07);

\draw [thick] (-60:2) .. controls (0,-1.1) and (-1,0) .. (-.5,.5);
\draw [thick] (1,0) .. controls (.5,.5) and (0,1) .. (-.5,.5);
\draw [thick] (1,0) .. controls (.5,.5) and (0,1) .. (-.5,.5);
\draw [thick] (-120:2) .. controls (-1.2,-.9) and (-1.4,-.7) .. (-1.4,0);
\draw [thick] (1,0) .. controls (.9,1.5) and (-1.4,1.6) .. (-1.4,0);
\draw [thick] (-120:2) .. controls (-1.6,-.9) and (-1.7,-.4) .. (-1.7,0);
\draw [thick] (-1.7,0) .. controls (-1.7,2.1) and (1.5,2) .. (1.5,0);
\draw [thick] (0,0) .. controls (.5,-.7) and (1.5,-1) .. (1.5,0);

\draw [purple!50!blue, thick] (-120:2) .. controls (-.8,-1) and (-1,-.3) .. (-1,0);
 \end{scope}
 
 \end{tikzpicture}
\caption{Surgery on a system of tethers using a fixed tether $t$}\label{fig:surgery}
\end{center}
\end{figure}

\begin{proof} We choose a single fixed tether $t$, then use a surgery flow to deform $T$ into the star of the vertex $t$. The flow will decrease the complexity of a system $\sigma$ (in normal form with respect to~$t$), which we define to be the total number of points in the intersection of the interiors of $\sigma$ and~$t$.  

If $s$ is a tether which intersects $t$ at an interior point, let $x$ be the intersection point which is closest along $t$ to the end $b_i$ of $t$.  Perform surgery on  $s$ by cutting it at $x$ and moving both new endpoints down to $b_i$ (see Figure~\ref{fig:surgery}).  This creates two new arcs which can be isotoped to be disjoint from $s$ except at their endpoints. 
One of these arcs joins $b_i$ to a puncture, and one joins $b_i$ to some (possibly different) $b_j$.  Define $\Delta s$ to be the arc connecting $b_i$ to a puncture. Note that $\Delta s$ has smaller complexity than $s$. 

The conditions of Lemma~\ref{lemma:flow} are now met, with $X=T$ and the star of $t$ as the subcomplex $Y$, by defining $v_\sigma$ to be the tether in $\sigma$ containing the point of $int(\sigma)\cap int(t)$ closest to $b_i$ along $t$. Thus $T$ deformation retracts to the star of $t$, which is contractible, hence $T$ is contractible.
\end{proof}

A system of tethers $ \tau =\{t_1,\cdots,t_k\}$ is {\it coconnected\/} if the complement  $D-\tau$ is connected.  Note that a system is coconnected if and only if each arc in the system ends at a different puncture. Let $T^0=T^0_{n,d}$ be the subcomplex of $T_{n,d}$ consisting of isotopy classes of coconnected tether systems.  

\begin{proposition}\label{doublebad} 
The complex $T^0=T^0_{n,1}$ is contractible.
\end{proposition}

\begin{proof}  We prove that $T^0$ is contractible by induction on the number $n$ of punctures.  If $n=1$ then $T^0$ is a single point.  For the induction step we will use a link argument  (Corollary~\ref{lemma:subcomplex}) to show the inclusion map $T^0\incl T=T_{n,1}$ is a homotopy equivalence, so we need to specify which simplices of $T$ are bad. We define a simplex of $T$  to be {\em bad\/} if each tethered puncture has at least two tethers.    We check (1) every simplex in $T$ which is not in $T^0$ has a bad face, and  (2) if $\sigma$  and $\tau$  are two bad faces of a simplex of $T$, the join $\sigma*\tau$ is also bad.   

If $\sigma$ is a bad simplex, we also need to identify the subcomplex $G_\sigma$ of $\lk(\sigma)$ consisting of simplices which are good for $\sigma$.  In our case $\tau\in \lk(\sigma)$ is good for $\sigma$ if and only if $\tau$ consists of single tethers to punctures which are not used by $\sigma$.  The subcomplex $G_\sigma$  decomposes as a join $G_\sigma=T^0(P_1)*T^0(P_2)*\cdots * T^0(P_r)$, where $P_1,\cdots,P_r$ are the components of the space obtained by cutting $D$ open along $\sigma$ and each $T^0(P_i)$ is either empty or isomorphic to $T^0_{n_i,d_i}$ for some $n_i<n$. Two tethers in $\sigma$ going to the same puncture bound a disk in $D$. A minimal such disk must be a component $P_i$ with at least one puncture in its interior (since isotopic tethers are not allowed) and only one distinguished boundary point.  Thus  $T^0(P_i)\cong T^0_{n_i,1}$ is contractible by induction on $n$, and the entire join $G_\sigma$ is contractible.  The hypotheses of  Corollary~\ref{lemma:subcomplex}(a) are satisfied and we conclude that $T^0$ is contractible since $T$ is.
\end{proof}

\begin{thm} \label{thm:braidstability}
The stabilization $H_i(B_{n-1})\to H_i(B_n)$ is an isomorphism for $n>2i+1$ and a surjection for $n=2i+1$. 
\end{thm}

\begin{proof} 
We use the spectral sequence constructed in Section~\ref{sec:spectral} for the action of $B_n$ on the contractible complex $T^0=T^0_{n,1}$.  Recall that this action arises from regarding $B_n$ as the group of isotopy classes of diffeomorphisms of the disk that are the identity on the boundary and permute the punctures $p_i$.  We verify  conditions (1)--(3) at the beginning of Section~\ref{sec:spectral}.

\begin{enumerate}
\item[(1)] $T^0$ has dimension $n-1$ and the action of $B_n$ has only one orbit of $k$-simplices for each $k$.  
\item[(2)] To see that the stabilizer of a $k$-simplex fixes the simplex pointwise, note that a set of $k+1$ tethers coming out of the basepoint in the boundary of the disk has a natural ordering determined by an orientation of the disk at the basepoint, and this ordering is preserved by any diffeomorphism of the disk that is the identity on the boundary. The stabilizer of a $k$-simplex is therefore isomorphic to $B_{n-k-1}$.  
\item[(3)] For an edge of $T^0$ corresponding to a pair of tethers there is a diffeomorphism of the disk supported in a neighborhood of the two tethers that interchanges the punctures at the ends of the tethers and takes the first tether to the second or vice-versa. This diffeomorphism gives an element of $B_n$ commuting with the stabilizer of the edge. 
\end{enumerate}
Part (a) of Theorem~\ref{thm:specseq} now 
gives the result since $T^0$ is contractible by Proposition~\ref{doublebad}.
\end{proof}
In this simple example we can in fact deduce more from the spectral sequence without much work:
\begin{thm} \label{thm:braidparttwo}
When $n$ is odd the stabilization $H_i(B_{n-1})\to H_i(B_n)$ is an isomorphism for all $i$. Also, $H_i(B_n)=0$ for $i\ge n$ (for $n$ of either parity).
\end{thm}
\begin{proof}

We look more closely at the spectral sequence used in the proof of Theorem~\ref{thm:braidstability} above.  The $E^1$ page has the following form:

\[
\begin{tikzpicture}
   \matrix (m) [matrix of math nodes,
                      nodes in empty cells, row sep=1em, column sep=4ex, minimum width=10ex]
  {
 H_i(B_n) &   H_i(B_{n-1}) &  H_i(B_{n-2}) & \cdots \hspace{3.5em}\\ 
 \hspace{1em} \cdots\hspace{.5em} &  H_{i-1}(B_{n-1}) &    H_{i-1}(B_{n-2})   & H_{i-1}(B_{n-3})&\cdots \hspace{3.5em}  \\ 
 &  &   \hspace{3.5em}{\cdots} \hspace{.5em}   & H_{i-2}(B_{n-3}) &\cdots \hspace{3.5em} \\ [3em]
  H_0(B_n) &   H_0(B_{n-1}) &  H_0(B_{n-2}) & H_0(B_{n-3})&\cdots\hspace{3.5em} \\ 
 p=\text{--}1 &  0 & 1 & 2&\cdots\hspace{3.5em} \\ 
};
  \path[-stealth]
    (m-1-2) edge   (m-1-1) 
    (m-1-3) edge node [above] {$0$}(m-1-2)
    (m-1-4) edge node [above] {$\iso$}(m-1-3)
   (m-2-2) edge   (m-2-1)
    (m-2-3) edge node [above] {$0$} (m-2-2)
    (m-2-4) edge node [above] {$\iso$}(m-2-3)
    (m-2-5) edge node [above] {$0$}(m-2-4)
    (m-3-4) edge node [above] {$\iso$} (m-3-3)
    (m-3-5) edge node [above] {$0$}(m-3-4)
    (m-4-2) edge (m-4-1)
    (m-4-3) edge node [above] {$0$} (m-4-2)
    (m-4-4) edge node [above] {$\iso$}  (m-4-3)
     (m-4-5) edge node [above] {$0$} (m-4-4);
\node [left] (i) at (-6.5,2.55) {$i$};
\node [left] (ione) at (-6.5,1.55) {$i{-}1$};
\node [left] (itwo) at (-6.5,.55) {$i{-}2$};
\node [left] (q) at (-6.5,-1.55) {$q = 0$};
\draw (-6.45,-3) -- (-6.45,3);
\draw (-7.5,-2.2) -- (5.25,-2.2);
\end{tikzpicture}
\]

To see that the differentials are alternately zeros and isomorphisms as shown, note first that the observation we used in the proof of Theorem~\ref{thm:braidstability} to verify condition (3) holds more generally to show that for any system $\sigma$ of  $k+1\geq 2$ tethers there is a diffeomorphism of the disk permuting the punctures and supported in a neighborhood of the tethers that takes any subset of $k$ of the tethers to any other set of $k$ of the tethers, preserving their natural order and commuting with $stab(\sigma)$.  This implies that each of the $p+1$ terms of the map $d^1\colon E^1_{p,q}=H_q(stab(\sigma))\to  E^1_{p-1,q}$ is the same, so for $p$ odd $d^1$ is zero, and for $p$ even $d^1$ is the map induced by inclusion.
If we assume $n$ is odd then by induction on $n$ starting with the trivial case $n=1$ the differential $d^1\colon E^1_{p,q} \to  E^1_{p-1,q}$  is an isomorphism for $p$ even, $p>0$. In particular, for the right-most nonvanishing column, which is the $p=n-1$ column since $T^0$ has dimension $n-1$, the $d^1$ differentials originating in this column are isomorphisms since $n$ is odd. (The only nonzero term in this column is $H_0(B_0)=\Z$ since $B_0$, like $B_1$, is the trivial group.) 

Thus in the $E^2$ page all the terms to the right of the $p=0$ column vanish.  Since the spectral sequence converges to zero and no differentials beyond the $E^1$ page can be nonzero, it follows that the differentials $d^1\colon H_i(B_{n-1})\to H_i(B_n)$ must be isomorphisms for all $i$, which finishes the induction step to prove the first part of the theorem.  

For the second statement of the theorem we again look at the $E^1$ page of the spectral sequence. The groups along the diagonal $p+q=n-1$ are the groups $H_j(B_j)$.  By induction on $n$ all the terms on or above this diagonal are zero except possibly in the $p=-1$ column, where the groups on or above the diagonal are $H_i(B_n)$ for $i\geq n$.  Since the spectral sequence converges to zero, all of these terms must vanish as well.
\end{proof}

The fact that $H_i(B_n)$ vanishes for $i\ge n$ is also a consequence of the well-known fact that there is an Eilenberg--MacLane space $K(B_n,1)$ which is a CW complex of dimension $n-1$ (see, e.g.,  \cite{CD}). 
Arnold proved the statements in the two preceding theorems in \cite{A} by methods not involving spectral sequences. 

Arnold also computed the homology of the pure braid subgroup  $P_n \subset B_n$ in \cite{A2} and it does not stabilize, even at the level of $H_1$ which is free abelian of rank $n\choose 2$, as can be seen already from a presentation for $P_n$. When the action of $B_n$ on $T^0$ is restricted to $P_n$ it is no longer transitive on simplices of each dimension, and in particular not on vertices.  We note that homology stability can sometimes still be proved using an action which is not transitive on simplices, as long as the number of orbits is independent of the stabilization parameters.  However, the spectral sequence argument becomes more complicated if there is more than one orbit.

\section{Curve and arc complexes}\label{sec:curvecomplexes}

For a compact orientable surface $S=S_{g,s}$ of genus $g$ with $s$ boundary components the classical {\it curve complex\/} $C(S)$ has as its vertices the isotopy classes of embedded curves (circles) in $S$ which are nontrivial, i.e., do not bound a disk and are not isotopic to a component of $ \bdry S$. A set of vertices of $C(S)$ spans a simplex if the corresponding curves can be isotoped to be all disjoint, so they form a {\it curve system}.  We will be particularly interested in the subcomplex $C^0(S)$ whose simplices are the isotopy classes of coconnected curve systems, i.e. systems with connected complement.  We will show that $C^0(S)$ is highly connected by showing that $C(S)$ is highly connected and using a link argument to deduce the connectivity for $C^0(S)$.  These results are due originally to Harer in~\cite{Har1} and we follow the same overall strategy while simplifying the proofs of several of the individual steps. 

\subsection{Curves on surfaces with nonempty boundary}

To prove $C(S)$ is highly connected when $\bdry S$ is nonempty the idea is to compare $C(S)$ with three other complexes in a sequence
\[
A(S,\bdry_0 S)\supset A_\infty(S,\bdry_0 S)
\simeq \Sub(S,\bdry_0 S) \simeq C(S)  
\]
The case that $S$ is closed will be deduced from the non-closed case.

We start by defining $A(S,\bdry_0 S)$. An {\it arc system\/} on a bounded surface $S$ is a set of disjoint embedded arcs with endpoints on the boundary $\bdry S$, such that no arc is isotopic to an arc in $\bdry S$ and no two arcs in a system are isotopic to each other, where all isotopies of arcs are required to keep their endpoints in $\bdry S$.  We choose a component $\bdry_0 S$ of $ \bdry S$ and define the complex $A(S,\bdry_0 S)$  as the geometric complex whose $k$-simplices are the isotopy classes of systems of $k+1$ arcs whose endpoints all lie in $\bdry_0 S$.

\begin{proposition} \label{prop:fullarccomplex}
The  complex  $A(S,\bdry_0 S)$ is contractible whenever it is nonempty, i.e., when $S$ is not a disk or annulus.
\end{proposition}

\begin{proof}  This is an application of Lemma~\ref{lemma:flow}, using surgery to flow into the star of a fixed ``target" arc $a$. The complexity of a system that intersects $a$ minimally within its isotopy class is defined as the number of intersection points with $a$. To do the surgery we first choose an orientation for $a$.  An arc $b$ crossing $a$ is cut into two arcs at the point where it meets $a$ nearest the terminal point of $a$, and the two new endpoints are moved to this terminal point to produce a new arc system $\Delta b$ meeting $a$ in one fewer point than $b$.  The function $\sigma\mapsto b_\sigma$ assigns to a system $\sigma$ the arc of $\sigma$ meeting $a$ at the point closest to the terminal point of $a$.  
\end{proof}

We define an arc system in $A(S,\bdry_0S)$ to be {\em at infinity}  if it has some complementary component which is neither a disk nor an annular neighborhood of a boundary component.  (The terminology comes from the fact, observed by Harer, that arc systems at infinity can be identified with rational points in the boundary of the Teichm\"uller space of the surface.)    Arc systems at infinity form a subcomplex $A_\infty(S,\bdry_0S)$.  A  calculation using Euler characteristics shows that it takes at least $(2g+s-1)$ arcs to cut $S$ into disks and annuli when $s=1$, so in this case $A_\infty(S,\bdry_0S)$ contains the entire $(2g+s-3)$-skeleton of $A(S,\bdry_0 S)$. When $s>1$ there are similar statements with $2g+s-1$ replaced by $2g+s-2$ and $2g+s-3$ replaced by $2g+s-4$. 
The inclusion $A_\infty(S,\bdry_0S)\incl A(S,\bdry_0S)$ induces an injection on $\pi_i$ when $A_\infty(S,\bdry_0S)$ contains the $(i+1)$-skeleton of $A(S,\bdry_0S)$, so we deduce:

\begin{corollary} \label{cor:Ainfinity}
 $A_\infty(S,\bdry_0 S)$ is $(2g+s-4)$-connected if  $s=1$ and $(2g+s-5)$-connected if  $s>1$.
\end{corollary}

Now define the {\it subsurface complex\/} $\Sub(S,\bdry_0 S)$ to be the geometric realization of the poset $\widehat\Sub(S,\bdry_0 S)$ of isotopy classes of compact connected subsurfaces $F$ of $S$ such that one component of 
$\bdry F$ is $\bdry_0 S$ and the other components of $\bdry F $ that are not contained in $\bdry S$ form a nonempty curve system in $S$, possibly containing parallel copies of the same curve. In particular, no component of $\bdry F-\bdry S$ bounds a disk in $S$ or is isotopic to a component of $\bdry S$. 

To each arc system $\alpha$ with $ \bdry\alpha \subset \bdry_0 S$ we
can associate a subsurface $F(\alpha)$ of $S$ by first taking a regular neighborhood $N$ of $ \alpha \cup \bdry_0 S$ and then adjoining any components of $S-N$ that are disks or annuli with one boundary circle contained in $\bdry S$  
(see Figure~\ref{fig:arcs}).
\begin{figure}
\begin{center}
\begin{tikzpicture}[yscale=.7, xscale=.95]
\draw (-1,1.5) to (3,1.5);  
\draw (-1,-1.5) to (3,-1.5);
\fill [blue!20] (-1,0) ellipse (.5cm and 1.5cm);
\fill [blue!20](-1,1.5) to (3,1.5) to (3,-1.5) to (-1,-1.5); 
\draw [fill=blue!10] (3,0) ellipse (.5cm and 1.5cm);
\node (bdry) at (3.8,-1) {$\partial_0S$};
\begin{scope}[xshift=1.5cm]
\fill  [white] (62:.5) arc (62:-62:.5);
\draw   (75:.5) arc (75:-75:.5);
\draw [fill=white] (60:.5) arc (120:240:.5);;
\end{scope}
\begin{scope}[xshift=-.5cm]
\fill  [white] (62:.5) arc (62:-62:.5);
\draw   (75:.5) arc (75:-75:.5);
\draw [fill=white] (60:.5) arc (120:240:.5);;
\end{scope}
\draw  [thick] (2.6,1) arc (80:280:1.75 and 1); 
 \draw  [densely dotted] (3.5,.2) arc (45:180:.9 and 0.2); 
 \draw  (3.5,.2) arc (45:110:.9 and 0.2); 
 \draw [thick] (2,0.1) arc (180:245):.9 and 0.2);
   
\draw [thick,->] (4.5,0) to (5.5,0);

\begin{scope}[xshift=8cm]
\draw (-1,-1.5) to (3,-1.5);
\draw (-1,1.5) to (3,1.5);  
\fill [blue!20] (-1,0) ellipse (.5cm and 1.5cm);
\fill [blue!20](-1,1.5) to (3,1.5) to (3,-1.5) to (-1,-1.5); 
\draw [fill=blue!10] (3,0) ellipse (.5cm and 1.5cm);
\node (bdry) at (3.8,-1) {$\partial_0S$};
\draw [ultra thin, fill=red!10] (2.7, 1.5) arc (90:-90:.5cm and 1.5cm) to (3,-1.5)  arc (-90:90:.5cm and 1.5cm) to (2.7,1.5);
\fill [red!10] (2.45,.01) arc (110:60:.9 and 0.2) to  
   (3.3,0) to (3.3,.3) arc (45:110:.9 and 0.2); 
\draw  [ultra thin] (2.45,.01) arc (110:61:.9 and 0.2); 
\draw  [ultra thin] (2.45,.35) arc (110:62:.9 and 0.2); 

\draw  [ultra thin]  (2.6,.9) arc (80:280:1.60 and .9);
\draw [ultra thin]  (2.6,-1.2)  
     arc (-80:-280:1.85 and 1.2); 
 \draw [ultra thin, fill=red!40] (2.7,1.5) arc (90:270:.5cm and 1.5cm) to (3, -1.5)  arc (-90:-270:.5cm and 1.5cm) to (2.7,1.5);
 \fill[red!40]  (2.6,.9) arc (80:280:1.60 and .9) to (2.6, -1.3) to (2.6,-1.2)   
     arc (-80:-280:1.85 and 1.2) to (2.6,.9); 
\fill [red!40] (1.95,0.2) arc (180:225):.9 and 0.2) to (2.21,-.29) arc (-135:-180:.9 and 0.2);
  \draw [ultra thin] (1.95,-.15) arc (180:225):.9 and 0.2);
 \draw [ultra thin] (1.95,0.2) arc (180:225):.9 and 0.2);
  \draw [ultra thin] (1.95,-.15) arc (180:225):.9 and 0.2);
 \draw [thin]  (3,0) ellipse (.5cm and 1.5cm);
 \draw  [densely dotted] (1.89,.05) arc (-180:-227:.9 and 0.3); 
\begin{scope}[xshift=1.5cm]
\fill  [white] (62:.5) arc (62:-62:.5);
\draw   (75:.5) arc (75:-75:.5);
\draw [fill=white] (60:.5) arc (120:240:.5);;
\end{scope}
\begin{scope}[xshift=-.5cm]
\fill  [white] (62:.5) arc (62:-62:.5);
\draw   (75:.5) arc (75:-75:.5);
\draw [fill=white] (60:.5) arc (120:240:.5);;
\end{scope}
\end{scope}
\end{tikzpicture}\caption{The map $\widehat A_\infty(S,\bdry_0 S)\to \widehat\Sub(S,\bdry_0 S)$}\label{fig:arcs}
\end{center}
\end{figure}
Thus the simplices of $A_\infty(S,\bdry_0 S)$ correspond to systems $\alpha$ for which $F(\alpha) \ne S$, and $\alpha\mapsto F(\alpha)$ is a map   $f\colon \widehat A_\infty(S,\bdry_0 S)\to \widehat\Sub(S,\bdry_0 S)$ where $ \widehat A_\infty(S,\bdry_0 S)$ denotes the poset of simplices in $A_\infty(S,\bdry_0 S)$. 
This map is a poset map since   $\alpha\subset\beta$ implies $F(\alpha)\subset F(\beta)$.

\begin{proposition}  \label{prop:subsurfaces}
The map $f\colon \widehat A_\infty(S,\bdry_0 S)\to \widehat\Sub(S,\bdry_0 S)$ is a homotopy equivalence.
\end{proposition}

\begin{proof}
We apply Quillen's fiber lemma, Proposition~\ref{lemma:poset}.  If $F$ is a subsurface of $S$ then  $f_{\leq F}$  is all arc systems $\alpha$ with $F(\alpha)\subset F$, so this is  $\widehat A(F,\bdry_0 S)$.  Since $F$ is not a disk or annulus,  $\widehat A(F,\bdry_0 S)$ is contractible by Proposition~\ref{prop:fullarccomplex}. 
\end{proof}

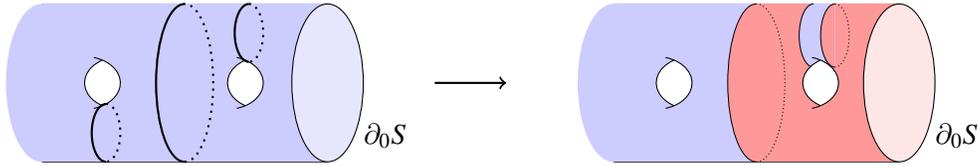
\begin{figure}
\begin{center}
\begin{tikzpicture}[yscale=.7, xscale=.95]
\draw (-1,1.5) to (3,1.5);  
\draw (-1,-1.5) to (3,-1.5);
\fill [blue!20] (-1,0) ellipse (.5cm and 1.5cm);
\fill [blue!20](-1,1.5) to (3,1.5) to (3,-1.5) to (-1,-1.5); 
\draw [fill=blue!10] (3,0) ellipse (.5cm and 1.5cm);
\node (bdry) at (3.8,-1) {$\partial_0S$};
\begin{scope}[xshift=1.6cm]
\fill  [white] (62:.5) arc (62:-62:.5);
\draw   (75:.5) arc (75:-75:.5);
\draw [fill=white] (60:.5) arc (120:240:.5);;
\end{scope}
\begin{scope}[xshift=-.4cm]
\fill  [white] (62:.5) arc (62:-62:.5);
\draw   (75:.5) arc (75:-75:.5);
\draw [fill=white] (60:.5) arc (120:240:.5);;
\end{scope}
\draw  [thick] (1,1.5) arc (90:270:.4 and 1.5);
\draw [thick, dotted] (1,-1.5) arc (-90:90:.4 and 1.5);
\draw  [thick] (1.9,1.5) arc (90:270:.2 and .55);
\draw [thick, dotted] (1.9,1.5) arc (90:-90:.2 and .55);
\draw   [thick] (-.1,-.4) arc (90:270:.2 and .55);
\draw [thick, dotted] (-.1,-.4) arc (90:-90:.2 and .55);

\draw [thick,->] (4.5,0) to (5.5,0);

\begin{scope}[xshift=8cm]
\draw (-1,-1.5) to (3,-1.5);
\draw (-1,1.5) to (3,1.5);  
\fill [blue!20] (-1,0) ellipse (.5cm and 1.5cm);
\fill [blue!20](-1,1.5) to (1,1.5) to (1,-1.5) to (-1,-1.5); 
\node (bdry) at (3.8,-1) {$\partial_0S$};
\fill [red!40] (1,0) ellipse (.4cm and 1.5cm);
\fill [red!40] (1,1.5) to (3,1.5) to (3,-1.5) to (1,-1.5);
\draw [fill=red!10] (3,0) ellipse (.5cm and 1.5cm);
\draw  [ultra thin] (1,1.5) arc (90:270:.4 and 1.5);
\draw [densely dotted] (1,-1.5) arc (-90:90:.4 and 1.5);
\begin{scope}[xshift=-.4cm] 
\fill  [white] (62:.5) arc (62:-62:.5);
\draw   (75:.5) arc (75:-75:.5);
\draw [fill=white] (60:.5) arc (120:240:.5);;
\end{scope}
\draw  [ultra thin, fill=blue!20] (1.8,1.5) arc (90:270:.2 and .6);
\draw [densely dotted, fill=blue!20] (1.8,1.5) arc (90:-90:.2 and .6);
\fill [blue!20] (1.8,1.5) to (2.1,1.5) to (2.1,.3) to (1.8,.3);
\draw  [fill=red!40, ultra thin] (2.1,1.5) arc (90:270:.2 and .6);
\draw [densely dotted] (2.1,1.5) arc (90:-90:.2 and .6);
\begin{scope}[xshift=1.65cm]
\fill  [white] (62:.5) arc (62:-62:.5);
\draw   (75:.5) arc (75:-75:.5);
\draw [fill=white] (60:.5) arc (120:240:.5);;
\end{scope}
\end{scope}
\end{tikzpicture}
\caption{The map $\widehat C(S) \to \widehat\Sub(S,\bdry_0 S)$}\label{fig:curves}
\end{center}
\end{figure}

Given a curve system $\gamma$, let the subsurface  $F(\gamma) \subset S$ be the component of the complement of a regular neighborhood of $\gamma$ containing $\bdry_0 S$ (see Figure~\ref{fig:curves}). 
Note that if $\gamma \subset \gamma'$ then $ F(\gamma)\supset F(\gamma')$. Thus if $\widehat C(S)$ denotes the poset of simplices of $C(S)$, then the association $ \gamma \mapsto F(\gamma)$ defines a poset map $g\colon \widehat C(S) \to \widehat\Sub(S,\bdry_0 S)$ with respect to the reverse ordering on $\widehat\Sub(S,\bdry_0 S)$ defined by $F_1\leq F_2$ if $F_1 \supset F_2$.

\begin{proposition} \label{prop:curvesandsubsurfaces}
The map $g\colon \widehat C(S) \to \widehat\Sub(S,\bdry_0 S)$ is a homotopy equivalence.
\end{proposition}

\begin{proof} We again apply Proposition~\ref{lemma:poset}.  For a subsurface $F$ in $\widehat\Sub(S,\bdry_0 S)$, the fiber $g_{\leq F}$ consists of curve systems in $S-F$, where curves are allowed to be parallel to curves of the system $\gamma(F) = \bdry F - \bdry S$.  In particular, $\gamma(F)$ is in the fiber, and  $\gamma(F)$ can be added to any curve system in the fiber, so the poset maps $\gamma\mapsto\gamma\cup \gamma(F)\mapsto \gamma(F)$ give a deformation retraction of $g_{\leq F}$ to the point $\gamma(F)$. 
\end{proof}

\begin{corollary} \label{cor:curvecomplex}
If $\bdry S$ is not empty, then $C(S)$ is $(2g+s-4)$-connected if $s=1$ and $(2g+s-5)$-connected if $s>1$.
\end{corollary}

\begin{corollary} \label{cor:genuszero}
If $S$ has genus $0$, then $C(S)$ is homotopy equivalent to a wedge of spheres of dimension $s-4$.
\end{corollary}

\begin{proof} If $S$ has genus $0$,  $C(S)$ is $(s-5)$-connected by the preceding corollary, and it is $(s-4)$-dimensional, so it is homotopy equivalent to a wedge of spheres of dimension $s-4$.
\end{proof}

\noindent{\bf Remark.}  In fact $C(S)$ is homotopy equivalent to a wedge of spheres in all cases.  When $g>0$ the dimension of the spheres is $2g+s-3$ if $s>0$ and $2g-2$ if $s=0$.  
This was proved by Harer in~\cite{Har2}, Theorems~3.3 and 3.5.  Thus the connectivity statements derived above are best possible when $s=1$ but one below best possible when $s>1$. 
However, these stronger results are not needed for the proof of homology stability.  

\subsection{Curves on closed surfaces}

There is a map $\phi\colon C(S_{g,1}) \to C(S_{g,0})$ induced by filling in the boundary circle of $S_{g,1}$ with a disk.  We remark that the dimension of $C(S_{g,1})$ is one more than that of $C(S_{g,0})$ when $ g > 1$ since maximal curve systems cut $S$ into pairs of pants. For $g=1$ the map $C(S_{1,1})\to C(S_{1,0})$ is an isomorphism.

\begin{proposition} \label{thm:curvecomplexclosedsurface}
The map  $\phi\colon C(S_{g,1})\to C(S_{g,0})$ is a homotopy equivalence for each $g\geq 1$.  
\end{proposition}

The weaker statement that $\phi_*\colon \pi_k C(S_{g,1})\to \pi_k C(S_{g,0})$ is surjective for all $k$ suffices to prove that $C(S_{g,0})$ is $(2g-3)$-connected, which is all we will need for homology stability.  The surjectivity of $\phi_*$ has a short proof using a little hyperbolic geometry, as follows.  Choose a hyperbolic structure on $S_{g,0}$ in the nontrivial cases $g\geq 2$.  Given a map $f\colon S^k\to C(S_{g,0})$ which we may assume is simplicial in some triangulation of $S^k$, the images $f(v)$ of all the vertices $v$ in $S^k$ can be represented by geodesics.  These are unique in their homotopy classes and are disjoint for sets of vertices spanning a simplex in $f(S^k)$.  Then a lift of $f$ to $C(S_{g,1})$ is obtained by deleting a disk in $S_{g,0}$ disjoint from this finite set of geodesics.

\smallskip
To obtain the full strength of Proposition~\ref{thm:curvecomplexclosedsurface} here is a proof that uses only topological techniques: 

\begin{proof}
We may assume $g\geq 2$.  It will suffice to show that for each simplex $\sigma$ of $C(S_{g,0})$ the subcomplex $F_\sigma = \phi^{-1}(\sigma)$ of $C(S_{g,1})$ is contractible, by Proposition~\ref{lemma:poset} and Lemma~\ref{lemma:simplicialmap}.  To begin, choose a curve system $\widetilde \sigma$ in $S_{g,1}$ with $\phi(\widetilde\sigma)=\sigma$. Enlarge $\widetilde\sigma$ to a maximal curve system $\delta$ cutting $S_{g,1}$ into pairs of pants. Let $P$ be the pair of pants containing $\bdry S_{g,1}$ and let $d_1$ and $d_2$ be the other two boundary circles of $P$. We can choose $\delta$ so that $d_1$ is a curve of $\widetilde\sigma$.

We may assume that all curve systems $\gamma$ in $S_{g,1}$ are in normal form with respect to $\delta$, so $\gamma$ intersects $\delta$ transversely in the minimum number of points among all systems isotopic to $\gamma$.  This minimality is equivalent to the ``no bigon'' condition that $S$ contains no disk whose boundary consists of an arc in $\gamma$ and an arc in $\delta$.
If two systems in normal form with respect to $\delta$ are isotopic, then they are isotopic through systems transverse to $\delta$, except that curves in $\gamma$ isotopic to curves in $\delta$ can be pushed from one side of $\delta$ to the other and such an isotopy cannot be transverse to $\delta$ at all times. 

If a curve system $\gamma$ is in normal form with respect to $\delta$ then each component arc of $\gamma \cap P$  either crosses $ P$ from $d_1$ to $d_2$, or it enters $P$, goes around $\bdry S_{g,1}$, and leaves by crossing the same $d_i$ that it crossed when it entered $P$.  An arc of the latter type we call a {\it return arc}. Note that all return arcs of $\gamma$ must have their endpoints on the same~$d_i$.

We will use surgery to flow from $F_\sigma$ into the subcomplex $F_\sigma^{nr}$ consisting of curve systems with no return arcs. 
Let $c$ be a curve in normal form with respect to $\delta$ that contains return arcs. Let $b$ be the innermost of these return arcs, the one closest to $\bdry S_{g,1}$. Pushing $b$ across $\bdry S_{g,1} $ converts $c$ into a new curve $\Delta c$ which can be isotoped to be disjoint from $c$. (See Figure~\ref{fig:stepone}). Alternatively, we can view $\Delta c$ as the result of surgering $c$ along an arc of $\bdry P$ to produce two curves, one of which is isotopic to $\bdry S_{g,1}$ and is discarded.  The curve $\Delta c$ may not be in normal form with respect to $\delta$, but it can be made so by an isotopy eliminating bigons one by one.  Note that $\Delta c$ is in $F_\sigma$ if $c$ is since the two curves become isotopic when $\bdry S_0$ is capped off with a disk.

If $\gamma$ is a curve system in $F_\sigma$ with at least one return arc, define $c_\gamma$ to be the curve in $\gamma$ containing the innermost return arc of $\gamma$.  Pushing $c_\gamma$ across $\bdry S_{g,1}$ as above then yields the curve $\Delta c_\gamma$.  If we define the complexity of a curve system to be the number of intersection points with $\bdry P$, we then have the ingredients to apply Lemma~\ref{lemma:flow}, producing a deformation retraction of $F_\sigma$ to $F_\sigma^{nr}$.  
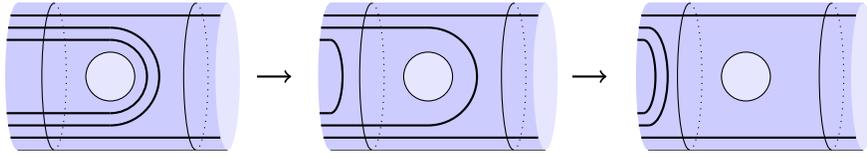
\begin{figure}
\begin{center}
\begin{tikzpicture}[scale=.65]  
\draw (-1.9,1.5) to (2.4,1.5);  
\draw (-1.9,-1.5) to (2.4,-1.5);
\fill [blue!20] (-1.9,0) ellipse (.25cm and 1.5cm);
\fill [blue!20](-1.9,1.5) to (2.4,1.5) to (2.4,-1.5) to (-1.9,-1.5); 
\fill [blue!10] (2.4,0) ellipse (.25cm and 1.5cm);
\draw   (1.75,1.5) arc (90:270:.25cm and 1.5cm);
\draw  [dotted] (1.75,1.5) arc (90:-90:.25cm and 1.5cm);
\draw   (-1.15,1.5) arc (90:270:.25cm and 1.5cm);
\draw  [dotted] (-1.15,1.5) arc (90:-90:.25cm and 1.5cm);
 \draw [fill=blue!10] (0,0) circle (.5cm);
 \draw [thick] (0,.75) arc (90:-90:.75 and .75);
 \draw [thick](-2.12,.75) to (0,.75);
\draw [thick](-2.12,.-.75) to (0,-.75);
  \draw [thick](0,1) arc (90:-90: 1 and 1);
   \draw [thick](-2.09,1) to (0,1);
\draw [thick](-2.09,.-1) to (0,-1);
\draw [thick] (-2.04,1.25) to (2.25, 1.25);
\draw [thick] (-2.04,-1.25) to (2.25, -1.25);

\draw [thick , ->] (3.0,0) to (3.7,0); 
  \begin{scope}[xshift=6.5cm]  
\draw (-2,1.5) to (2.4,1.5);  
\draw (-2,-1.5) to (2.4,-1.5);
\fill [blue!20] (-2,0) ellipse (.25cm and 1.5cm);
\fill [blue!20](-2,1.5) to (2.4,1.5) to (2.4,-1.5) to (-2,-1.5); 
\fill [blue!10] (2.4,0) ellipse (.25cm and 1.5cm);
\draw   (1.75,1.5) arc (90:270:.25cm and 1.5cm);
\draw  [dotted] (1.75,1.5) arc (90:-90:.25cm and 1.5cm);
\draw   (-1.15,1.5) arc (90:270:.25cm and 1.5cm);
\draw  [dotted] (-1.15,1.5) arc (90:-90:.25cm and 1.5cm);
 \draw [fill=blue!10] (0,0) circle (.5cm);
 \draw [thick] (-2,.75) arc (90:-90:.25 and .75);
 \draw [thick](-2.22,.75) to (-2,.75);
\draw [thick](-2.22,.-.75) to (-2,-.75);
  \draw [thick](0,1) arc (90:-90: 1 and 1);
   \draw [thick](-2.19,1) to (0,1);
\draw [thick](-2.19,.-1) to (0,-1);
\draw [thick] (-2.14,1.25) to (2.25, 1.25);
\draw [thick] (-2.14,-1.25) to (2.25, -1.25);
\draw [thick , ->] (2.95,0) to (3.65,0); 
 \end{scope}
 
  \begin{scope}[xshift=13cm] 
\draw (-2,1.5) to (2.4,1.5);  
\draw (-2,-1.5) to (2.4,-1.5);
\fill [blue!20] (-2,0) ellipse (.25cm and 1.5cm);
\fill [blue!20](-2,1.5) to (2.4,1.5) to (2.4,-1.5) to (-2,-1.5); 
\fill [blue!10] (2.4,0) ellipse (.25cm and 1.5cm);
\draw   (1.75,1.5) arc (90:270:.25cm and 1.5cm);
\draw  [dotted] (1.75,1.5) arc (90:-90:.25cm and 1.5cm);
\draw   (-1.15,1.5) arc (90:270:.25cm and 1.5cm);
\draw  [dotted] (-1.15,1.5) arc (90:-90:.25cm and 1.5cm);
 \draw [fill=blue!10] (0,0) circle (.5cm);
 \draw [thick] (-2.1,.75) arc (90:-90:.25 and .75);
 \draw [thick](-2.22,.75) to (-2.1,.75);
\draw [thick](-2.22,.-.75) to (-2.1,-.75);
  \draw [thick](-2,1) arc (90:-90: .4 and 1);
   \draw [thick](-2.19,1) to (-2,1);
\draw [thick](-2.19,.-1) to (-2,-1);
\draw [thick] (-2.14,1.25) to (2.25, 1.25);
\draw [thick] (-2.14,-1.25) to (2.25, -1.25);
\end{scope}

 \end{tikzpicture}

\caption{The flow in the proof of Proposition~\ref{thm:curvecomplexclosedsurface}}\label{fig:stepone}
 \end{center}
 \end{figure}

We claim that $F_\sigma^{nr}$ is just a single simplex, the simplex spanned by $\widetilde \sigma$ and $d_2$. To see this, observe that if $\gamma$ is a simplex in $F_\sigma^{nr}$ in normal form with respect to $\delta$, and with any of its curves parallel to $d_1$ or $d_2$ pushed outside $P$, then we can obtain the normal form for $\phi(\gamma)$ by simply deleting $P$ from $S_{g,1}$ and identifying $d_1$ and $d_2$ in such a way as to match up the endpoints of any arcs of $\gamma \cap P$.  In fact there can be no such arcs since $\phi(\gamma)$ is a face of $\sigma$, hence $\gamma$ is a face of $\widetilde\sigma*d_2$.  Thus $F_\sigma^{nr}$ is a simplex and it follows that $F_\sigma$ is contractible.
\end{proof}

\subsection{Coconnected curve systems}

The complex $C^0(S)$ of coconnected curve systems on a surface of genus $g$ has dimension ${g-1}$.  As we will observe in Remark~\ref{rem:tophomology} below, the top-dimensional homology group $H_{g-1}(C^0(S))$ is nonzero, so the best one could hope is that $C^0(S)$ is $(g-2)$-connected, and indeed it is:

\begin{proposition}  \label{thm:coconnectedcurvesystems}
The complex  $C^0(S)$ of coconnected curve systems on a surface $S$ of genus $g$ is $(g-2)$-connected.
\end{proposition}

\begin{proof}   
This is a link argument, an application of Corollary~\ref{lemma:subcomplex} with $X=C(S)$ and $Y=C^0(S)$. To begin we need to single out the bad simplices of $C(S)$. 
To each curve system we associate a dual graph, with a vertex for each complementary component of the system and an edge for each curve.  Thus a curve system is coconnected if and only if its dual graph has one vertex and all edges are loops.  A {\em bad simplex\/} in $C(S)$ is a system of curves for which {\em no} edges of the dual graph are loops.  This is equivalent to saying that each curve in the system separates the complement of the other curves.  It is easy to see that conditions (1) and (2) in Section~\ref{subsec:linksarguments} are satisfied for this notion of badness.  For a bad simplex $\sigma$, the complex $G_\sigma$ is the join of the complexes $C^0(S_i)$ for the components $S_i$ of the surface $S_\sigma$ obtained by cutting $S$ open along $\sigma$. Either the genus $g_i$ of $S_i$ is smaller than $g$ or $S_i$ has fewer boundary components than $S$ so we may proceed by induction on the lexicographically ordered pair $(g,s)$.  
Since $g_\sigma=\sum_ig_i$ by definition and the quantity ``connectivity plus two" is additive for joins, it follows that we may assume inductively that $G_\sigma$ is $(g_\sigma -2)$-connected. The induction can start with the obvious cases that the genus is zero or one.

Cutting $S$ along the curves of $\sigma$ can decrease the total genus by at most $\dim(\sigma)$ since each cut decreases genus by at most one and cutting along the last curve cannot decrease genus since $\sigma$ is bad.  Thus $g_\sigma\geq g-\dim(\sigma)$.  (This estimate is best possible in the case that the dual graph to $\sigma$ consists of two vertices joined by a number of edges.)  
It follows that the connectivity of $G_\sigma$, which inductively is at least $g_\sigma -2$, is also at least $g-2-\dim(\sigma)$ so the connectivity hypothesis on links in Corollary~\ref{lemma:subcomplex}(a) is satisfied with $n=g-2$.

The last thing to check to apply the corollary is that the larger complex $C(S)$ is $(g-2)$-connected.  We can assume $g>0$ since the proposition is trivial when $g=0$.  Then Corollary~\ref{cor:curvecomplex} and Proposition~\ref{thm:curvecomplexclosedsurface} imply that $C(S)$ is $(2g-3)$-connected, and we have $g-2\leq 2g-3$ when $g\geq 1$.
\end{proof}

\begin{remark}\label{rem:tophomology}{\rm 
 There is an easy argument showing that $H_{g-1}(C^0(S_{g,s}))$ is nonzero for all $g\geq 1$ and $s\geq 0$.  Choosing $g$ disjoint copies of $S_{1,1}$ in $S_{g,s}$ gives an embedding of the join of $g$ copies of $C^0(S_{1,1})$ into $C^0(S_{g,s})$ as a subcomplex. The complex $C^0(S_{1,1})=C(S_{1,1})$ is an infinite discrete set, so the join is homotopy equivalent to the wedge of an infinite number of copies of $S^{g-1}$.  The inclusion map of the join into $C^0(S_{g,s})$ induces an injection on $H_{g-1}$ since both complexes have dimension $g-1$ so no nontrivial $(g-1)$-dimensional cycle in the join can bound in $C^0(S_{g,s})$. Thus $H_{g-1}(C^0(S_{g,s}))$ is nontrivial, and in fact is free abelian of infinite rank since it is the kernel of the boundary map from the free abelian group of simplicial $(g-1)$-chains to the simplicial $(g-2)$-chains, and a subgroup of a free abelian group is free abelian.}
 \end{remark}

There is an oriented version of $C^0(S)$ whose simplices are isotopy classes of coconnected systems of curves together with choices of orientations for these curves.  Call the resulting complex $C^0_\pm(S)$.

\begin{corollary} \label{cor:orientedCzero}
The complex $C^0_\pm(S)$ is $(g-2)$-connected.
\end{corollary}

\begin{proof}
Choose an arbitrary orientation for each isotopy class of nonseparating curves in $S$. Then the two possible orientations correspond to the labels $+$ and $-$ and the result is immediate from Corollary~\ref{coloring}.
\end{proof}

\begin{remark} {\rm There are also versions of $C^0(S)$ and $C^0_\pm(S)$ in which simplices correspond to {\em ordered\/} coconnected systems of curves or oriented curves.  These too have the same connectivity as $C^0(S)$ by Proposition~\ref{prop:OrderedComplex}. 
 }
\end{remark}

\section{ Tethered curves and chains}
\label{sec:tethers}

This section represents the heart of the paper, where we introduce the geometric complexes that encode more information than is given by curves or arcs alone. The main work is in showing that the new complexes are roughly half as highly-connected as $C^0(S)$, but this is enough for the spectral sequence arguments.
The various complexes we will consider fit into a commutative diagram:
\[
\begin{tikzpicture}
\matrix (m) [matrix of math nodes, row sep=1.5em, column sep=2.2em, minimum width=2em]
{TCh(S) & DTC(S) & TC(S) \\
Ch(S) &&  C^0(S) \\};
\path[-stealth]
   (m-1-1) edge (m-1-2) edge  (m-2-1)
   (m-1-2) edge (m-1-3)
   (m-1-3) edge  (m-2-3)
   (m-2-1) edge  (m-2-3);
\end{tikzpicture}
\]
\noindent
The maps are forgetful maps except for the upper left horizontal map which is an injection.  We will start with the known connectivity of the complex  $C^0(S)$ in the lower right corner, then proceed around the diagram in the counterclockwise direction to show each complex in turn is highly connected.  Except for the one injection, each step will involve two stages: first enlarge the domain complex to a complex for which a surgery flow can be used to show that the fibers of the extended forgetful map are contractible, then use a link argument to shrink back to the original source complex.  For the injection we need only a link argument to show that the image (and hence the domain) is highly connected.

\subsection{Tethered curves}
\label{subsec:tethers}

Let $P$ be a nonempty finite collection of disjoint open intervals and circles in $\bdry S$. 
A {\it tether\/} for a simple closed curve $c$ in $S$ is an arc in $S$ with one endpoint in $c$ and the other in $P$, the interior of the arc being disjoint from $c$ and from $\bdry S$.  Define a complex $TC(S,P)$ whose $k$-simplices are isotopy classes of systems of $k+1$ disjoint tethered curves such that the complement of the system of tethered curves is connected.  This last condition is equivalent to the curves by themselves forming a coconnected system since after cutting $S$ open along the curves, each newly created boundary circle is connected to $P$ by at most one tether arc, so cutting along these arcs cannot disconnect the surface.  Note that tethering a curve gives it a normal orientation, pointing away from curve in the direction of the tether. An orientation of $S$ converts the normal orientation of the curve to a tangential orientation. 

 \begin{proposition}\label{prop:TCconnectivity}
The complex $TC(S,P)$ is $(g-3)/2$-connected, where $g$ is the genus of $S$.
\end{proposition}
\begin{proof} 
We will include $TC(S,P)$ into a larger complex $TC(S,P^\ast)$, then project this larger complex  onto $C^0_\pm(S)$, the oriented version of $C^0(S)$.  We show that the projection is a homotopy equivalence using a surgery argument on barycentric fibers. This shows that $TC(S,P^\ast)$ is $(g-2)$-connected, and then a link argument will show that $TC(S,P)$ is $(g-3)/2$-connected. 

The complex $TC(S,P^\ast)$ has the same vertices as $TC(S,P)$ but the simplices correspond to coconnected curve systems with at least one but possibly more tethers to each curve.  All the tethers must be disjoint, and all tethers to the same curve must attach on the same side of the curve (at distinct points), giving the curve the same orientation.   

The projection $TC(S,P^\ast)\to C^0_\pm(S)$ forgets the tethers but remembers the orientations they induce on curves.
 The fiber of this projection over the barycenter of a simplex $\sigma$ of $C^0_\pm(S)$ consists of all tether systems for $\sigma$ attaching on the ``positive" sides of the curves in $\sigma$. Choose one such system, which for simplicity we take to lie in $TC(S,P)$ so that it consists of a single tether $t_i$ to each curve $c_i$ in $\sigma$. 
 We can deform the barycentric fiber into the star of this tethered system by a surgery flow.  The surgeries are performed first using the tether $t_1$, surgering toward $P$ until all tethers are disjoint from $t_1$, then using $t_2$ in similar fashion, and so on.  Each surgery cuts a tether into two arcs, one of which has both ends in $P$ and which we discard, while the other is a tether which we keep.  To make the surgery process well defined on isotopy classes one must first put the tethers being surgered into normal form with respect to the fixed tethers $t_i$; this minimizes the number of intersection points with $t_i$ by eliminating bigons and ``half-bigons". 
 
The surgery flow shows that the barycentric fiber over $\sigma$ is contractible, so the projection $TC(S,P^\ast)\to C^0_\pm(S)$ is a homotopy equivalence by Proposition~\ref{lemma:fiber}.  Thus $TC(S,P^\ast)$ is $(g-2)$-connected by Corollary~\ref{cor:orientedCzero}.
 
We now use a link argument to analyze the inclusion of $TC(S,P)$ into $TC(S,P^\ast)$.  
To see how the connectivity number $(g-3)/2$ arises, let us try to show the inclusion is $(Ag+B)$-connected for yet-to-be-determined constants $A$ and $B$.  A {\em bad simplex\/} in $TC(S,P^\ast)$ is one that corresponds to a system of tethered curves in which each curve has at least two tethers; in particular a vertex cannot be bad.  If $\sigma$ is a bad simplex, the surface $S_\sigma$ obtained by cutting $S$ open along the system of curves and tethers corresponding to $\sigma$ may have several components $S_i$, and $P$ is cut into pieces $P_i\subset \bdry S_i$.  The components $S_i$ all have smaller genus than $S$ and the subcomplex  $G_\sigma$  is the join of the complexes $TC(S_i,P_i)$,  so we can argue by induction on the genus. 

To apply Corollary~\ref{lemma:subcomplex}(a) we need to arrange that $G_\sigma$ is $(Ag+B - k)$-connected for $k=\dim(\sigma)$.  Cutting a surface along a multi-tethered curve system decreases genus by at most one for each tether, so if $g_\sigma$ denotes the genus of $S_\sigma$ (i.e. the sum of the genera $g_i$ of the components $S_i$) we have $g_\sigma \geq g-k-1$. 
Suppose we know by induction on genus that $G_\sigma$ is $(Ag_\sigma+B)$-connected.  The inequality we need is $Ag_\sigma+B\geq Ag+B-k$, i.e., $Ag_\sigma \geq Ag-k$.  We have observed that $g_\sigma\geq g-k-1$, so it suffices to have $A(g-k-1)\geq Ag-k$, which simplifies to $A\leq k/(k+1)$.  We only need this for $k\geq 1$ since vertices cannot be bad, so $A=1/2$ works for all $k\geq1$.  We therefore choose $A=1/2$. 

To apply Corollary~\ref{lemma:subcomplex}(a) we also need $TC(S,P^\ast)$ to be $(g/2 +B)$-connected, which means $g/2 +B\leq g-2$, the connectivity of $TC(S,P^\ast)$. This inequality reduces to $B\leq g/2-2$. We can assume $g\geq 1$ since the proposition is trivially true when $g=0$.  When $g=1$ the inequality $B\leq g/2-2$ is $B\leq -3/2$ so we maximize $B$ by choosing $B=-3/2$ and then $B\leq g/2-2$ for all $g\geq 1$.

Thus our candidate for $Ag+B$ is $(g-3)/2$.  It remains to verify the induction step by showing that the join of the $(g_i-3)/2$-connected complexes $TC(S_i, P_i)$ is $(g_\sigma -3)/2$-connected.  
We know that connectivity plus two is additive for joins, but this assumes the connectivities are integers and here they could be fractions.  This means we need to use the floor function $\lfloor \cdot \rfloor$ for connectivities in order to apply the connectivity-plus-two fact. Thus we let $$f(g)=\bigl\lfloor  (g-3)/2\bigr\rfloor+2$$ and we wish to verify that $f(g_\sigma)\leq \sum_i f(g_i)$. 

 We have $f(g)=g/2+1/2$ if $g$ is odd and $g/2$ if $g$ is even.  If $g_\sigma=\sum_ig_i$ is odd then at least one $g_i$ is odd, so $\sum_i f(g_i)\geq \bigl(\sum_i g_i/2\bigr)+1/2 = f(g_\sigma)$.  If $g_\sigma$ is even, we only need to notice that $\sum_i f(g_i)\geq  \sum_i g_i/2=f(g_\sigma)$.
\end{proof}

\subsection{Double-tethered curves}
\label{subsec:doubletethers}

\begin{figure}
\begin{center}
\begin{tikzpicture}[yscale=.7, xscale=.95]

 \draw[densely dotted]  (0,0)  [partial ellipse=-90:90:.25cm and 1cm]; 
 \draw  (0,0)  [partial ellipse=90:270:.25cm and 1cm]; 
 \draw[densely dotted]  (0,3)  [partial ellipse=-90:90:.25cm and 1cm]; 
 \draw  (0,3)  [partial ellipse=90:270:.25cm and 1cm]; 
\draw  (5,3) ellipse (.25cm and 1cm); 
\draw (5,0) ellipse (.25cm and 1cm); 

\draw [very thick, blue] (5,0) [partial ellipse=150:210:.25cm and 1cm]; 
\draw [very thick, blue] (5,3) [partial ellipse=150:210:.25cm and 1cm];
 
\draw (0,1) arc ( -90 : 90 : .5); 
\draw (5,2) arc ( 90 : 270 : .5); 
\draw (0,-1)  .. controls (1,-.5) and (4,-.5) ..  (5,-1); 
\draw (0,4)  .. controls (1,3.5) and (4,3.5) ..  (5,4); 

 \coordinate (x1) at (4.75,0);
 \coordinate (x2) at (4.75,3);
 \node [right] (xtwo) at (x2) {$P$}; 
 \fill (x2) circle (.05);
  \node [right] (xone) at (x1) {$Q$};
   \fill (x1) circle (.05);
   
\begin{scope}[xshift=2.5cm, yshift=.9cm]  
 \fill  [white] (-28:.75) arc (-28:-152:.75);
 \draw   (-15:.75) arc (-15:-165:.75);
 \draw [fill=white] (-30:.75) arc (30:150:.75);
 \end{scope}
\begin{scope}[xshift=2.5cm, yshift=2.9cm] 
 \fill  [white] (-28:.75) arc (-28:-152:.75);
 \draw   (-15:.75) arc (-15:-165:.75);
 \draw [fill=white] (-30:.75) arc (30:150:.75);
 \end{scope}

\draw[thick, red, densely dotted] (2.5, 1.52) [partial ellipse=-90:90:.2cm and .62cm];
\draw[thick, red] (2.5, 1.52) [partial ellipse=90:270:.2cm and .62cm];
\draw[thick, red] (x2) .. controls (0,4.5) and (1,1.5) .. (2.2,1.62);  
\draw[thick, red] (2.2,1.62).. controls (4,1.7) and (4,0) .. (x1);
\end{tikzpicture}
\caption{Double-tethered curve}\label{fig:doubletether} 
\end{center}
\end{figure}
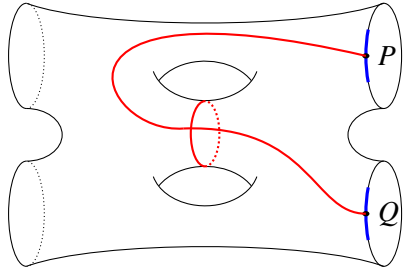

 We now consider a complex $DTC(S,P,Q)$ of double-tethered curve systems.  Here $Q$ is a second nonempty finite collection of disjoint open intervals and circles in $\bdry S$ (we allow $P$ and $Q$ to overlap or even coincide), and a {\em double tether} for a curve $c$ is an ordered pair of tethers attaching at the same point of $c$ but on opposite sides, with the first tether going to a point in $P$ and the second to a point in $Q$.  The two tethers must be disjoint except at their common attaching point in $c$ (see Figure~\ref{fig:doubletether}).  It is often useful to think of a double tether as a single oriented arc going from $P$ to $Q$ and crossing the curve at a single point.  Note that a Dehn twist along the curve $c$ acts nontrivially on the isotopy classes of double tethers for $c$, in contrast with the situation for single tethers. 
 
A $k$-simplex of $DTC(S,P,Q)$ is by definition an isotopy class of systems of $k+1$ disjoint double-tethered curves such that the complement of the system is connected.  As before, this last condition is equivalent to the curves by themselves forming a coconnected system since cutting along the tethers after cutting along the curves cannot disconnect the surface. 
 
\begin{proposition}\label{prop:DTCconnectivity}
The complex $DTC(S,P,Q)$ is $(g-3)/2$-connected, where $g$ is the genus of $S$.
\end{proposition}

\begin{proof}  
The proof follows closely the proof of Proposition~\ref{prop:TCconnectivity}.  As before we include $DTC(S,P,Q)$ into a larger complex $DTC(S,P,Q^\ast)$, then project to $TC(S,P)$.
The complex $DTC(S,P,Q^\ast)$ has the same vertices as $DTC(S,P,Q)$; however higher-dimensional simplices of $DTC(S,P,Q^\ast)$ correspond to isotopy classes of coconnected curve systems with one tether from each curve to $P$ and at least one but possibly more tethers to $Q$, where all the tethers for a given curve attach at the same point of the curve and all the $Q$-tethers attach on the opposite side from the $P$-tether.  All the tethers to $P$ and $Q$ in the system are disjoint from each other and from the curves except at the points where they attach to a curve.  Faces of simplices in $DTC(S,P,Q^\ast)$ are obtained by deleting one $Q$-tether to a curve if there are several, or by deleting the whole double-tethered curve if there is only one $Q$-tether to it. The vertices of the simplex are thus the double-tethered curves contained in the given system of curves and tethers.

The projection $DTC(S,P,Q^\ast)\to TC(S,P)$ forgets the tethers to $Q$, keeping only the single tether from each curve to $P$.  As before, the barycentric fiber over a simplex of $TC(S,P)$ can be contracted by surgery into the star of a fixed system in $DTC(S,P,Q)$.  Thus $DTC(S,P,Q^\ast)$ is $(g-3)/2$-connected.

We now use a link argument exactly as in the proof of Proposition~\ref{prop:TCconnectivity} to deduce that $DTC(S,P,Q)$  
is $(g-3)/2$-connected. 
The key point is that cutting $S$ along a simplex of $DTC(S,P,Q^\ast)$  decreases the genus by at most one for each tether to $Q$.  This is because cutting along a nonseparating curve decreases genus by one, then cutting along the single tether to $P$ does not decrease the genus further nor does cutting along the first tether to $Q$, and cutting along each additional tether to $Q$ can decrease genus by at most one.
\end{proof}

For the spectral sequence proof of homology stability we will use a certain subcomplex of $DTC(S,P,Q)$ defined when $P$ and~$Q$ are disjoint single intervals.  To define this subcomplex we first choose orientations for $P$ and~$Q$.  For a simplex of $DTC(S,P,Q)$ the orientation of $P$ induces an ordering of the double tethers of this simplex. Likewise the orientation of $Q$ induces a possibly different ordering of the double tethers.  The simplices for which the two orderings are in fact the same form a subcomplex of $DTC(S,P,Q)$ which we denote $DTC^{m}(S,P,Q)$, with the superscript indicating {\it matching\/} orderings. 

\begin{proposition}\label{prop:ODTCconnectivity}
The complex $DTC^{m}(S,P,Q)$ is $(g-3)/2$-connected.
\end{proposition}

\begin{proof} This will be a link argument, following the idea of the proof of Theorem~4.9 of \cite{W2}.  A simplex of $DTC(S,P,Q)$ has vertices a set of pairs $(c_0,d_0),\cdots,(c_k,d_k)$ consisting of curves $c_i$ with double tethers $d_i$.  We may assume these are listed in the order specified by the orientation of $P$.  The ordering determined by the orientation of $Q$ differs from this ordering by a permutation $\pi$ of $\{0,1,\cdots,k\}$, with $\pi$ the identity exactly when the simplex is in $DTC^{m}(S,P,Q)$.  If $\pi$ is not the identity let $i$ be the smallest index such that $\pi(i)\neq i$.  We call the given simplex {\em bad\/} if $i=0$.  We can write the simplex uniquely as the join of a simplex $\langle(c_0,d_0),\cdots,(c_{i-1},d_{i-1})\rangle$ which is good, i.e., in $DTC^{m}(S,P,Q)$, and a simplex $\langle(c_i,d_i),\cdots,(c_k,d_k)\rangle$ which is bad, where either of these two subsimplices could be empty.  This notion of badness satisfies the two conditions in Subsection~\ref{subsec:linksarguments}. 

For a bad simplex $\sigma = \langle(c_0,d_0),\cdots,(c_k,d_k)\rangle$ the subcomplex $G_\sigma$ of simplices that are good for $\sigma$ can be identified with $DTC^{m}(S_\sigma,P_\sigma,Q_\sigma)$ where $S_\sigma$ is the (connected) surface obtained by cutting $S$ along $\sigma$ and $P_\sigma$ and $Q_\sigma$ are the subintervals of $P$ and $Q$ up to the point where the first (with respect to the orientations of $P$ and $Q$) tethers of $\sigma$ attach.  Cutting $S$ open along each double tethered curve decreases genus by one, so by induction on genus we may assume that $DTC^{m}(S_\sigma,P_\sigma,Q_\sigma)$ is $(g-(k+1)-3)/2$-connected. Since there are no bad $0$-simplices we have $k\geq 1$ and $(g-(k+1)-3)/2\geq (g-3)/2 - k$.  
The result now follows from Corollary~\ref{lemma:subcomplex}(a).
\end{proof}

\subsection{ Chains and tethered chains}\label{subsec:chains}

The connectivity results obtained so far are enough to prove homology stability for mapping class groups of surfaces with nonempty boundary.  However if a surface is closed there is no natural place for tethers to go, and we instead consider complexes of {\em chains}, where a   chain is an ordered pair $(a,b)$ of simple closed curves intersecting transversely in a single point, together with an orientation 
 on $b$.  The geometric complex of chains is denoted $Ch(S)$.

\begin{remark} {\rm Forgetting the orientation on $b$ gives a retraction of $Ch(S)$ to a complex of unoriented chains, which has the same connectivity by an application of Corollary~\ref{coloring}.}
\end{remark}

We will prove that $Ch(S)$ is highly-connected for all surfaces $S$, with or without boundary, but we start with a complex of tethered chains $TCh(S_{g,s},P)$ for a surface with non-empty boundary.  Here each chain has one tether connecting the positive side of the (oriented) $b$-curve of the chain to a point in some finite collection $P$ of disjoint open intervals (but no circles) in $\bdry S$. Since the tether is only allowed to attach to the positive side of the $b$-curve, specifying the tether determines the orientation of $b$. 

\begin{proposition}\label{prop:TChconnectivity}
The complex $TCh(S_{g,s},P)$ is $(g-3)/2$-connected.
\end{proposition}

\begin{proof}  
Let $((a,b),t)$ be a vertex of $TCh(S_{g,s},P)$.  A small neighborhood $N$ of $b\cup t$ is homeomorphic to an annulus, one of whose boundary components intersects $\bdry S$ in an arc.  Deleting the interior of this arc from this component of $\bdry N$ leaves a double tether $t'$ for the $a$-curve (see Figure~\ref{fig:btether}). %
\begin{figure}
\begin{center}
  \begin{tikzpicture}
  \fill[white] (0,-2) rectangle    (6,2);  
   \draw (0,1)  .. controls (1.5,.75) and (1.5,1.5) ..  (3,1.5); 
    \draw (3,1.5)  .. controls (4.5,1.5) and (4.5,.75) ..  (6,1); 
    Ê \draw (0,-1)  .. controls (1.5,-.75) and (1.5,-1.5) ..  (3,-1.5); 
     \draw (3,-1.5)  .. controls (4.5,-1.5) and (4.5,-.75) ..  (6,-1); 
      \fill [red!40, opacity=0.2] (3,0) ellipse (1.7cm and 1cm);
  \fill[white]  (3,0) ellipse (1.3cm and .6cm);
  \fill [red!40,opacity=0.2] (.253,-.2) to (1.35,-.2) to (1.35,.2) to (.253,.2);
    \draw  [thick] (.25,.2) to (1.35,.2);
\begin{scope}[xshift=3cm,yshift=.4cm]   
 \fill  [white] (-28:.75) arc (-28:-152:.75);
 \draw   (-15:.75) arc (-15:-165:.75);
 \draw [fill=white] (-30:.75) arc (30:150:.75);
 \end{scope}
 \draw [red, thick] (3,0) ellipse (1.5cm and .8cm);
 \draw[thick, red] (3,-.925) [partial ellipse=-90:90:.24cm and .575cm];
\draw[thick, red, densely dotted] (3,-.925) [partial ellipse=90:270:.24cm and .575cm];
 \draw [thick] (3,0) ellipse (1.7cm and 1cm);
  \draw  (3,0) ellipse (1.3cm and .6cm);
 \fill [white] (1.2,-.19) rectangle (1.4,.19);
\fill [red!40,opacity=.2] (1.2,-.19) rectangle (1.4,.19);
  \draw  [thick] (.25,-.2) to (1.35,-.2);
    \draw  [thick] (.25,.2) to (1.35,.2);

 \draw[thick, red] (.25,0) to (1.5,0);
 
  \draw  (0,0) ellipse (.25cm and 1cm);
 
\node (a) at (3.38,-1.28) {$a$};
 \fill [white, opacity=.75]  (4.5,0) ellipse (.1cm and .2cm);
\node (b) at (4.5,0) {$b$};
 \fill [white, opacity=.75]  (.8,0) ellipse (.1cm and .15cm);
\node (t) at (.8,0) {$t$};
\node (t') at (1.41,.68) {$t'$};
 \end{tikzpicture}
 \caption{Double tether associated to a tethered chain }
 \label{fig:btether}
 \end{center}
 \end{figure}
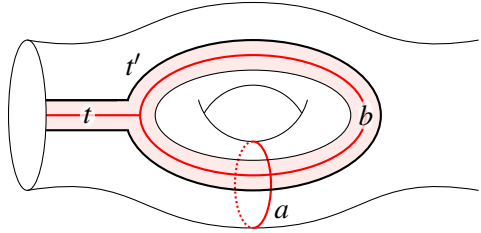
 
Thus the double tether $t'$ and the $a$-curve give a vertex of the complex $DTC(S_{g,s},P)$ of double-tethered coconnected curve systems where the double tethers have both ends in $P$.  In $DTC(S_{g,s},P)$ we do not orient the double tethers, in contrast with the double tethers in $DTC(S_{g,s},P,P)$ which do have specified orientations. 
 This map on vertices $((a,b),t)\mapsto (a,t')$ extends to a simplicial embedding $$TCh(S_{g,s},P)\hookrightarrow DTC(S_{g,s},P).$$  
The image consists of simplices with the special property that the two ends of each double tether used in the simplex are adjacent in $P$. We denote this image by $DTC^{a}(S_{g,s},P)$, with the superscript indicating the adjacency of the two ends of a double tether. 
Thus it suffices to prove that $DTC^{a}(S_{g,s},P)$ is $(g-3)/2$-connected.  We will do this by a link argument similar to the one for Proposition~\ref{prop:ODTCconnectivity}.  This will use the fact that the larger complex $DTC(S_{g,s},P)$ is $(g-3)/2$-connected, which follows by embedding it in $DTC(S_{g,s},P,P)$ as a retract by arbitrarily choosing orientations for all the double tethers of vertices of $DTC(S_{g,s},P)$; the retraction is obtained by replacing all orientations by these arbitrarily chosen ones, and $DTC(S_{g,s},P,P)$ is $(g-3)/2$-connected by Proposition~\ref{prop:DTCconnectivity}. 

For the link argument assume first that $P$ is a single interval and choose an orientation for~$P$.  This allows us to order the ends of the double tethers in each simplex of $DTC(S_{g,s},P)$. Define a simplex of $DTC(S_{g,s},P)$ to be {\em bad\/} if its first double-tether end in $P$ is not immediately followed by the other end of this double tether.  (Note that vertices cannot be bad.)  Each simplex in $DTC(S_{g,s},P)$ is then the join of two of its faces, the first face consisting of a string (possibly empty) of adjacently double-tethered curves whose tether ends form an initial segment of the sequence of all the tether ends, and the second face a bad simplex whose tether ends form the rest of the sequence. The two conditions for badness in Subsection~\ref{subsec:linksarguments} are easily checked.

For a bad $k$-simplex $\sigma$ of $DTC(S_{g,s},P)$ the subcomplex $G_\sigma$ of simplices that are good for $\sigma$ can be identified with $DTC^{a}(S_\sigma,P_\sigma)$ where $S_\sigma$ is the surface obtained by cutting $S$ along $\sigma$ and $P_\sigma$ is the part of $P$ up to the first attaching point for the tethers of $\sigma$.  By induction on genus we may assume $G_\sigma$ is $(g-(k+1)-3)/2$-connected,  
hence  $((g-3)/2 - k)$-connected since $k\geq 1$.
 The result for the case that $P$ is a single interval then follows from Corollary~\ref{lemma:subcomplex}(a). 

Now we treat the case of a more general $P$ consisting of several disjoint intervals. Let $P_0$ be one of these intervals. We will apply a link argument for $TCh(S_{g,s},P)$ and its subcomplex $TCh(S_{g,s},P_0)$.  Define a $k$-simplex $\sigma$ of  $TCh(S_{g,s},P)$ to be {\em bad\/} if all of its tethers attach to points in $P-P_0$.  Clearly the two conditions for badness are satisfied, and $G_\sigma$ is $TCh(S_\sigma,P_0)$ for $S_\sigma$ the surface obtained by cutting $S_{g,s}$ along $\sigma$. We have shown that $TCh(S_\sigma,P_0)$ is $(g-k-1-3)/2$-connected, hence $((g-3)/2 -k-1)$-connected.  By Corollary~\ref{lemma:subcomplex}(b), since $TCh(S_{g,s},P_0)$ is $((g-3)/2)$-connected, so is $TCh(S_{g,s},P)$.
\end{proof}

Finally we consider the complex  $Ch(S_{g,s})$ of oriented chains.
 
\begin{proposition}\label{prop:Chconnectivity}
 $Ch(S_{g,s})$ is $(g-3)/2$-connected. 
\end{proposition}

\begin{proof}
We first treat the cases $s>0$. Consider the complex $TCh(S_{g,s},P)$ with $P$ a collection of disjoint open intervals in $\bdry S_{g,s}$. We enlarge $TCh(S_{g,s},P)$ to a complex $TCh(S_{g,s},P^\ast)$ by allowing multiple tethers to each chain, all attaching at the same point of the $b$-curve of the chain and on the same side of the curve, each tether being otherwise disjoint from all other tethers and chains. There is a projection $$TCh(S_{g,s},P^\ast)\to Ch(S_{g,s})$$ obtained by forgetting the tethers and orienting chains according to which side of the $b$-curves the tethers attach to.  The fibers of this projection are contractible by the usual surgery argument, so it suffices to show that $TCh(S_{g,s},P^\ast)$ is $(g-3)/2$-connected.  

We do this by a link argument. The {\em bad\/} $k$-simplices $\sigma$ in $TCh(S_{g,s},P^\ast)$ are those whose chains all have at least two tethers.  The complex $G_\sigma$ is the join of the complexes $TCh(S_i,P_i)$ where cutting $(S_{g,s},P)$ along $\sigma$ produces a pair $(S_\sigma,P_\sigma)$ with components $(S_i,P_i)$.  Cutting along a chain reduces genus by one and creates a new boundary circle, and then cutting along the first tether to the chain does not reduce the genus further, while cutting along each subsequent tether to the chain reduces genus by at most one more.  Thus the total genus of $S_\sigma$ is at least $g-k-1$.  It follows as in the last paragraph of the proof of Proposition~\ref{prop:TCconnectivity} that $G_\sigma$ is $(g-k-1-3)/2$-connected, hence $((g-3)/2 -k-1)$-connected. 
Since $TCh(S_{g,s},P)$ is $((g-3)/2)$-connected by Proposition~\ref{prop:TChconnectivity}, we can apply Corollary \ref{lemma:subcomplex}(b) to deduce that $TCh(S_{g,s},P^\ast)$ has this connectivity as well.  This proves the  proposition when $s>0$.

For the case $s=0$ we use hyperbolic geometry.  The cases $g\leq 1$ are trivial, so we can assume $g\geq 2$ and fix a hyperbolic structure on $S_{g,0}$.  Each nontrivial isotopy class of curves in $S_{g,0}$ contains a unique geodesic representative, and the geodesic representatives of two isotopy classes intersect the minimum number of times within the isotopy classes.  Furthermore, if two curves intersect minimally, then one can choose isotopies to geodesics such that the number of intersections remains minimal throughout the isotopies.  Thus each simplex in $Ch(S_{g,0})$ has a unique geodesic representative.

There is a simplicial map $Ch(S_{g,1})\to Ch(S_{g,0})$ induced by filling in $\bdry S_{g,1}$ with a disk.  Given a simplicial map $f\colon S^i\to Ch(S_{s,0})$ we can choose a disk $D$ in $S_{g,0}$ disjoint from the finitely many geodesic representatives for the chains that are images of vertices of $S^i$.  Deleting the interior of $D$, we then have a lift of $f$ to $Ch(S_{g,1})$.  This lift is nullhomotopic if $i\leq (g-3)/2$.  Composing with the projection to $Ch(S_{g,0})$ then gives a nullhomotopy of $f$, so $Ch(S_{g,0})$ is $(g-3)/2$-connected.  
\end{proof}

\begin{remark} 
{\rm One may ask whether the case $s=0$ can be proved by a purely topological argument.   
In Proposition~\ref{thm:curvecomplexclosedsurface} this was done for the analogous projection $C(S_{g,1})\to C(S_{g,0})$ by showing (in essence) that its fibers, which are one-dimensional when $g\geq 2$, are contractible.  However, the fibers of $Ch(S_{g,1})\to Ch(S_{g,0})$ are zero-dimensional and infinite when $g\geq 2$ so we cannot expect the same approach to work here. }
\end{remark}

\begin{example}\label{ex:genustwochains} {\rm Consider the case that $S$ is closed of genus $2$, so $Ch(S)$ is one-dimensional.  A chain $(a,b)$ in $S$ has neighborhood bounded by a separating curve $c=c(a,b)$.  The link of $(a,b)$ in $Ch(S)$ consists of all chains $(a',b')$ in the genus one surface on the other side of $c$. These chains also have $c(a',b')=c$, so it follows that all chains $(a,b)$ in this connected component of $Ch(S)$ have the same curve $c(a,b)$.  The connected components of $Ch(S)$ thus correspond to nontrivial separating curves on $S$.  Each connected component is the join of two copies of the infinite zero-dimensional complex $Ch(S_{1,1})$. Thus $Ch(S)$ is not homotopy equivalent to a wedge of spheres of a single dimension, in contrast with the situation for the curve complexes $C(S)$ and $C^0(S)$.  Note also that the connectivity bound $(g-3)/2$ is best possible in this case, where it asserts only that $Ch(S)$ is nonempty.}  
\end{example}

\begin{remark}\label{rem:tophomologyagain} {\rm As a partial generalization of this example one can say that for $S$ a surface of arbitrary genus~$g\geq 1$ the group $H_{g-1}(Ch(S))$ is free abelian of infinite rank. This follows as in Remark~\ref{rem:tophomology} by embedding $g$ disjoint copies of $S_{1,1}$ in $S$, which gives an embedding of the join of $g$ copies of the infinite discrete set $Ch(S_{1,1})$ in $Ch(S)$.  This join has dimension $g-1$, the same dimension as $Ch(S)$, so the embedding of the join is injective on $H_{g-1}$. }
\end{remark}

\section{The stability theorems}\label{sec:MCGstability}

In the first part of this section we apply the spectral sequence for the action of the mapping class group of $S$ on a double-tethered curve complex to prove homology stability with respect to genus when the number of boundary components is fixed and nonzero.  As a bonus, the proof also shows that the stable homology groups do not depend on the number of boundary components of $S$ when this number is nonzero.

After this we show that the homology in the case of closed surfaces is isomorphic to that for non-closed surfaces, in the stable dimension range.  This uses a second spectral sequence, this one for the action of the mapping class group on the complex $\ord{Ch(S)}$ of ordered chains in $S$.  

\subsection{Surfaces with nonempty boundary}

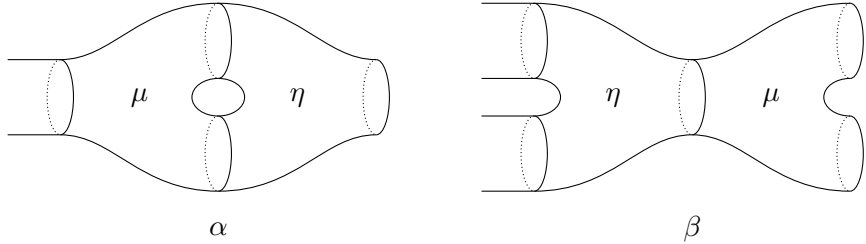
\begin{figure}
\begin{center}
\begin{tikzpicture}[xscale=.7, yscale=.5]
\begin{scope}[xshift=-12cm]

\draw (-4,1) to (-3,1);
\draw (-4,-1) to (-3,-1);
\draw (0,1.5) [partial ellipse=-90:90:.25cm and 1cm];
\draw [densely dotted](0,1.5) [partial ellipse=90:270:.25cm and 1cm];
\draw (0,-1.5) [partial ellipse=-90:90:.25cm and 1cm];
\draw [densely dotted](0,-1.5) [partial ellipse=90:270:.25cm and 1cm];
\draw (0,-.5) arc ( -90 : -270 : .5); 
 \draw (-3,0) [partial ellipse=-90:90:.25cm and 1cm];
 \draw[densely dotted] (-3,0) [partial ellipse=90:270:.25cm and 1cm];
\draw (0,2.5)  .. controls (-1.5,2.5) and (-2,1) ..  (-3,1); 
\draw (0,-2.5)  .. controls (-1.5,-2.5) and (-2,-1) ..  (-3,-1); 

 \node  (mu) at (-1.5,0) {$\mu$};
 
 \draw(3,0) [partial ellipse=-90:90:.25cm and 1cm];
 \draw [densely dotted] (3,0) [partial ellipse=90:270:.25cm and 1cm];
\draw (0,-.5) arc ( -90 : 90 : .5); 
\draw (0,2.5)  .. controls (1.5,2.5) and (2,1) ..  (3,1); 
\draw (0,-2.5)  .. controls (1.5,-2.5) and (2,-1) ..  (3,-1); 
 \node  (eta) at (1.5,0) {$\eta$};
 
  \node (alpha) at (0,-3.5) {$\alpha$};
 
 \end{scope}
 
 \begin{scope}[xshift=-6cm]
 
 \draw (-1,2.5) to (0,2.5);
  \draw (-1,.5) to (0,.5);
   \draw (-1,-.5) to (0,-.5);
    \draw (-1,-2.5) to (0,-2.5);
\draw (0,1.5) [partial ellipse=-90:90:.25cm and 1cm];
\draw [densely dotted](0,1.5) [partial ellipse=90:270:.25cm and 1cm];
\draw(0,-1.5) [partial ellipse=-90:90:.25cm and 1cm];
\draw [densely dotted] (0,-1.5) [partial ellipse=90:270:.25cm and 1cm];
\draw (0,-.5) arc ( -90 : 90 : .5); 
\draw (3,0) [partial ellipse=-90:90:.25cm and 1cm];
\draw [densely dotted](3,0) [partial ellipse=90:270:.25cm and 1cm];
\draw (0,2.5)  .. controls (1.5,2.5) and (2,1) ..  (3,1); 
\draw (0,-2.5)  .. controls (1.5,-2.5) and (2,-1) ..  (3,-1); 
 \node  (eta) at (1.5,0) {$\eta$};
 \end{scope}

\draw (0,1.5) [partial ellipse=-90:90:.25cm and 1cm];
\draw [densely dotted](0,1.5) [partial ellipse=90:270:.25cm and 1cm];
\draw (0,-1.5) [partial ellipse=-90:90:.25cm and 1cm];
\draw [densely dotted](0,-1.5) [partial ellipse=90:270:.25cm and 1cm];
\draw (0,-.5) arc ( -90 : -270 : .5); 
\draw (0,2.5)  .. controls (-1.5,2.5) and (-2,1) ..  (-3,1); 
\draw (0,-2.5)  .. controls (-1.5,-2.5) and (-2,-1) ..  (-3,-1); 
 \node  (mu) at (-1.5,0) {$\mu$};  

 \node (beta) at (-3,-3.5) {$\beta$};  
 
 \end{tikzpicture}
\caption{$\alpha$ and $\beta$ stabilizations}\label{fig:alphabeta}
\end{center}
\end{figure}

Let $S$ be a surface of genus $g$ with $s\ge 1$ boundary components, and let $M_{g,s}$ be the mapping class group of $S$, where diffeomorphisms and isotopies between them are required to restrict to the identity on each boundary circle. There are two stabilization maps $M_{g,s} \to M_{g+1,s}$ induced by inclusions $\alpha,\beta \colon S_{g,s} \to S_{g+1,s}$ shown in Figure~\ref{fig:alphabeta}.  For $\alpha$ one attaches $S_{1,2}$ to $S_{g.s}$ along one boundary circle, assuming $s\geq 1$, while for $\beta$ one attaches $S_{0,4}$ to $S_{g,s}$ along two boundary circles, assuming $s\geq 2$.  The inclusions $\alpha$ and $\beta$ induce homomorphisms of the corresponding mapping class groups by extending diffeomorphisms via the identity on the attached surface.  It is a standard fact that these induced homomorphisms are injective; see for example Theorem~3.18 in~\cite{FM}.  

We can factor both $\alpha$ and $\beta$ as compositions of two inclusions $\mu$ and $\eta$, each attaching a pair of pants $S_{0,3}$, with the $\mu$ attachment being along one boundary circle of $S_{0,3}$ and $\eta$ along two boundary circles.  The difference between $\alpha$ and $\beta$ is the order of the attachments: For $\alpha$ it is $\mu$ followed by $\eta$ while for $\beta$ it is the reverse.

For the associated mapping class groups we can stabilize with respect to $g$ for fixed $s$ by iterating $\alpha$ arbitrarily often, or we can do the same using $\beta$.  The $\alpha$ stabilization is the one usually considered rather than $\beta$, probably because it is the more obvious stabilization and only requires $s\geq 1$.  We can also iterate $\mu$ arbitrarily often to stabilize with respect to $s$ for fixed $g$, but $\eta$ can only be iterated finitely often so it is not exactly a stabilization.

\begin{thm} \label{thm:maintheorem}
The stabilizations $\begin{Bmatrix}
\alpha_*\colon H_i(M_{g,s})\to H_i(M_{g+1,s}) \\
\mu_*\colon H_i(M_{g,s})\to H_i(M_{g,s+1})
\end{Bmatrix}$ are isomorphisms for $g \geq 2i+2$ and $s\geq 1$.
\end{thm}

\begin{proof}  We first give an easy argument that  reduces both cases in the theorem to the statement that  
the stabilization $\beta_*\colon H_i(M_{g,s})\to H_i(M_{g+1,s})$ is an isomorphism for $g>2i+1$ and a surjection for $g=2i+1$.  Consider the following three maps:
$$
\begin{CD}
H_i(M_{g,s}) @>{\mu_*}>>  H_i(M_{g,s+1})  @>{\eta_*}>> H_i(M_{g+1,s})  @>{\mu_*}>> H_i(M_{g+1,s+1})
\end{CD}
$$
The composition of the first two maps is $\alpha_*$ and the composition of the second two is $\beta_*$. If $\beta_*$ is surjective for $g\geq 2i+1$ then so is the second $\mu_*$ in that range. On the other hand  $\mu_*$ is always injective since it has a left inverse obtained by attaching a disk to one of the free boundary circles of the attached $S_{0,3}$, so that the net result of the two attachments is attaching an annulus along one boundary circle, and this induces the identity map on $M_{g,s}$. 
Thus the second $\mu_*$ is an isomorphism for $g\geq 2i+1$.  It follows that $\eta_*$ is surjective for $g\geq 2i+1$.  The first $\mu_*$ is now an isomorphism for $g\geq 2i+2$, so $\alpha_*$ is surjective in that range.  For injectivity of $\alpha_*$, if $\beta_*$ is injective for $g\geq 2i+2$ then so is $\eta_*$, hence also $\alpha_*$ since $\mu_*$ is always injective.

To prove stability for $\beta_*$ we will apply part (b) of Theorem~\ref{thm:specseq} for the action of $M_{g,s}$ on the complex $DTC^{m}(S_{g,s},P,Q)$ of systems of double-tethered curves with matching orderings along $P$ and $Q$, where $P$ and $Q$ are single intervals in different components of $\bdry S_{g,s}$. Here we assume $s\geq 2$ in order for $\beta$ to be defined.  Since $P$ and $Q$ are single intervals we can use a slightly different, equivalent definition of $DTC^{m}(S_{g,s},P,Q)$ in which basepoints $x_1\in P$ and $x_2\in Q$ are chosen in advance and all tethers are required to have their $P$-endpoints at $x_1$ and their $Q$-endpoints at $x_2$, but otherwise satisfy the same conditions as before.  Note that orderings of the tethers at $x_1$ and $x_2$ are still well-defined, specified by orientations of $P$ and~$Q$.

We need to check that the conditions (1)--(3) at the beginning of Section~\ref{sec:spectral} hold for this action. First we check that the action is transitive on simplices of each dimension. The mapping class group clearly acts transitively on ordered coconnected systems of $k$ oriented curves. To see that this holds also when matched systems of double tethers are added, we use the orientation on a curve determined by its tethers which specify a $P$-side of the curve and a $Q$-side, and we use the ordering of the curves specified by the ordering of the tethers at $P$ and $Q$, which agree since we assume tethers satisfy the matching condition.  Then transitivity on the double tethers for a fixed coconnected ordered oriented curve system can be seen by first cutting $S$ along the curves in the system to get a surface $F$, and then observing that the mapping class group of $F$ acts transitively on systems of $k$ arcs in $F$ joining $x_1$ to basepoints $p_1,\cdots,p_k$ in the $k$ ordered $P$-circles of $\bdry F - \bdry S$ together with $k$ arcs joining $x_2$ to basepoints $q_1,\cdots,q_k$ in the $k$ ordered $Q$-circles of $\bdry F - \bdry S$, where in both cases the orderings of the arcs at $x_1$ and $x_2$ agree with the specified orderings of the circles at their other endpoints.  This can be seen inductively by first making any two arcs from $x_1$ to $p_1$ agree after a diffeomorphism, then making any two arcs from $x_1$ to $p_2$ starting on the same side of the first arc at $x_1$ agree after a diffeomorphism fixing the first arc, etc.
 
To see that the inclusion of the stabilizer of a vertex $\sigma$ into $M_{g,s}$ is the map induced by $\beta$, note first that a diffeomorphism in the stabilizer can be isotoped to fix the double-tethered curve pointwise, not just setwise, since it fixes $x_1$ and $x_2$.  Then it can be isotoped to be the identity in a closed neighborhood $N_\sigma$ of the union of $\sigma$ and the components of $\bdry S$ containing $x_1$ and $x_2$ (see Figure~\ref{fig:Nsigma}). This $N_\sigma$ is diffeomorphic to $S_{0,4}$ since it has Euler characteristic $-2$ and four boundary circles. Furthermore, $N_\sigma$ attaches to the complementary surface  $S_\sigma$  along two circles of $\bdry N_\sigma$. Thus the inclusion of the stabilizer of $\sigma$ is the $\beta$ stabilization. (If $x_1$ and $x_2$ were in the same component of $\bdry S$, the neighborhood $N_\sigma$ would be a copy of $S_{1,2}$ and the inclusion of the vertex stabilizer would be the $\alpha$ stabilization.)

  \begin{figure}
\begin{center}
\begin{tikzpicture}[yscale=.8]
\draw(-.25,0) [partial ellipse=-90:90:.25cm and 1cm];
\draw[densely dotted]  (-.25,0) [partial ellipse=90:270:.25cm and 1cm];
\draw[thick,densely dotted, red]  (0,0) [partial ellipse=120:240:.25cm and 1cm];
\draw(.25,0) [partial ellipse=-90:90:.25cm and 1cm];

\begin{scope}[xshift=-5cm]
\draw (0,0) [partial ellipse=-90:90:.25cm and 1cm];
\draw(0,0) [partial ellipse=90:270:.25cm and 1cm];
\draw(.25,0) [partial ellipse=-90:90:.25cm and 1cm];
\draw (.25,0) [partial ellipse=120:240:.25cm and 1cm];
\end{scope}

\begin{scope}[xshift=5cm]
\draw(-.25,0) [partial ellipse=-90:90:.25cm and 1cm];
\draw[densely dotted]  (-.25,0) [partial ellipse=90:270:.25cm and 1cm];
\draw (0,0) [partial ellipse=-90:90:.25cm and 1cm];
\draw[ densely dotted]  (0,0) [partial ellipse=120:240:.25cm and 1cm];
\end{scope}

\fill[color=white]   (-4.7,-.25) rectangle (5.2,.25);
\draw (-4.5,.25) to (0,.25);
 \draw (.5,.25) to (5,.25);
 \draw (-4.5,-.25) to (0,-.25);
  \draw (.5,-.25) to (5,-.25);
\draw [thick, red](-4.75,0) to (5.25,0);
\draw[thick, red](0,0) [partial ellipse=-90:90:.25cm and 1cm];

\draw (-5,1) ..controls (-4.5,1) and (-4.3,1).. (-4,1.5);
\draw (-5,-1) ..controls (-4.5,-1) and (-4.3,-1).. (-4,-1.5);

\draw (5,1) ..controls (4.5,1) and (4.3,1).. (4,1.5);
\draw (5,-1) ..controls (4.5,-1) and (4.3,-1).. (4,-1.5);

\draw (0,1) ..controls (.5,1) and (.7,1).. (1,1.5);
\draw (0,-1) ..controls (.5,-1) and (.7,-1).. (1,-1.5);

\draw (0,1) ..controls (-.5,1) and (-.7,1).. (-1,1.5);
\draw (0,-1) ..controls (-.5,-1) and (-.7,-1).. (-1,-1.5);

\fill [white, opacity=.85]  (-5,.03) ellipse (.1cm and .23cm);
\node (x1) [left] at (-4.63,0) {$x_1$};
\node (x2) [right] at (5.25,0) {$x_2$};
 \end{tikzpicture}
\caption{The neighborhood $N_\sigma$ for a vertex $\sigma$ of $DTC(S)$}\label{fig:Nsigma}
\end{center}
\end{figure}
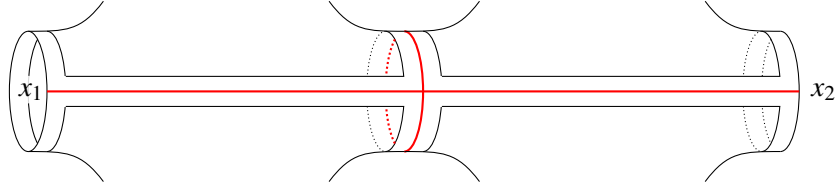

More generally the inclusion of the stabilizer of a $k$-simplex $\sigma$ is a $(k+1)$-fold iterate of $\beta$ stabilizations since cutting $S$ along each double-tethered curve of $\sigma$ in succession gives $k+1$ $\beta$ stabilizations.  All such inclusions of stabilizers of $k$-simplices into $M_{g,s}$ are conjugate since the action is transitive on $k$-simplices.

For condition (2), the stabilizer of a simplex fixes the simplex pointwise since the order of tethers at a basepoint cannot be permuted by an orientation-preserving diffeomorphism of the surface.  To check condition (3) note that a neighborhood of two double-tethered curves defining an edge of $DTC^{m}(S_{g,s})$, with the double tethers going from $x_1$ to $x_2$, is a copy of $S_{1,2}$ (this is not the same as the neighborhood $N_\sigma$ in the preceding paragraph since we do not include neighborhoods of the boundary circles containing $x_1$ and $x_2$). The mapping class group  $M_{1,2}$ acts transitively on vertices of $DTC^{m}(S_{1,2})$, so there is a diffeomorphism of $S$ supported in a neighborhood of the two given double-tethered curves that sends the first to the second, or vice versa.
\end{proof}

\begin{remark}  {\rm If we choose $P$ and $Q$ to be disjoint intervals in the same component of $\bdry S$, the same proof as above shows that the stabilization $\alpha_*\colon H_i(M_{g,s})\to H_i(M_{g+1,s})$ is an isomorphism for $g\geq 2i+2$ and a surjection for $g=2i+1$, a slight improvement over the preceding theorem.  
The advantage of using the $\beta$ stabilization is that one also gets $\mu$ stability for free.}
\end{remark}

\subsection{Closed surfaces}

It remains to deal with the projection $M_{g,1}\to M_{g,0}$ induced by filling in the boundary circle of $S_{g,1}$ with a disk.  More generally we will consider the map $\kappa\colon M_{g,s+1}\to M_{g,s}$ induced by capping off a boundary circle with a disk.

\begin{thm} \label{thm:closedsurfaces}
The map $\kappa_*\colon H_i(M_{g,s+1})\to H_i(M_{g,s})$ is an isomorphism for $g > 2i+3$ and a surjection for $g = 2i+3$, for each $s\geq 0$.  
\end{thm}

\begin{proof} 
We will apply part (b) of Theorem~\ref{thm:specseq} to the action of $M_{g,s}$ on $\ord{Ch(S_{g,s})}$, the complex of ordered systems of oriented chains.    
The orderings and orientations guarantee that the stabilizer of a $k$-simplex is exactly $M_{g-k-1,s+k+1}$.  The inclusion of this stabilizer is the $(k+1)$-fold iterate of the composition $\lambda=\kappa\alpha$.  We already know that $\alpha$ induces an isomorphism on homology in a stable range, so proving this for $\kappa$ is equivalent to proving it for $\lambda$ (and surjectivity of $\lambda_*$ implies surjectivity of $\kappa_*$). 

Conditions (1) and (2) for Theorem~\ref{thm:specseq} are obviously satisfied.  Unfortunately condition~(3) fails for $1$-simplices of $\ord{Ch(S_{g,s})}$ since there is no diffeomorphism of $S_{g,s}$ moving one chain onto another disjoint chain and supported in a neighborhood of the two chains.  However there is a weakening of condition~(3) that is satisfied and is strong enough to make the argument for injectivity of the differential $d\colon H_i(M_{g-1,s+1})\to H_i(M_{g,s})$ still work.  (The argument for surjectivity did not use (3).) 
If we enlarge the neighborhood of the two chains by adding a neighborhood of an arc joining them, producing a surface $T\subset S$ diffeomorphic to $S_{2,1}$, then there is a diffeomorphism supported in $T$ interchanging the two chains and preserving their orientations. 
If we denote the two vertices of $\ord{Ch(S_{g,s})}$ corresponding to the two chains by $v$ and $w$, with $e$ either of the two edges joining them, then we have the following diagram: 
\[
\begin{tikzpicture}
\matrix (m) [matrix of math nodes, row sep=1em, column sep=1.5em, minimum width=2em]
{ && H_{i}(stab(v)) &  \\
H_{i}(stab(T)) & H_{i}(stab(e)) & &H_{i}(stab(v_0)) \\
  && H_{i}(stab(w)) &  \\};
\path[-stealth]
   (m-2-1) edge [bend left=5] (m-1-3) edge  (m-2-2) edge [bend right=5] (m-3-3)
   (m-2-2) edge (m-1-3) edge   (m-3-3)
   (m-1-3) edge (m-2-4)  edge (m-3-3)
   (m-3-3) edge (m-2-4) 
  ;
\end{tikzpicture}
\]

\noindent
The four triangles in the diagram commute except possibly the one just to the left of the vertical arrow.
Also, the large triangle formed by the two curved arrows and the vertical map commutes. The horizontal map is the $\mu_*$ stabilization in Theorem~\ref{thm:maintheorem} so it is an isomorphism provided that $g-2\geq 2i+2$, i.e., $g\geq 2i+4$. 
This implies that the whole diagram is in fact commutative in this range.  
This suffices to deduce injectivity of the differential $d=\lambda \colon H_i(M_{g-1,s+1})\to H_i(M_{g,s})$ when $g> \varphi(i) = 2i+c$ for $c=3$, recalling that $c\geq 2$ was sufficient in the original argument for Theorem~\ref{thm:specseq}(b).  

It remains only to check that $\ord{Ch(S_{g,s})}$ is $(g-3)/2$-connected, which will follow from Proposition~\ref{prop:OrderedComplex} if $Ch(S_{g,s})$ is wCM of level $(g-1)/2$.  We know that $Ch(S_{g,s})$ is $(g-3)/2$-connected by Proposition~\ref{prop:Chconnectivity}, and likewise the link of a $k$-simplex of $Ch(S_{g,s})$ is $(g-(k+1)-3)/2$-connected. We have $(g-(k+1)-3)/2 \geq (g-1)/2 - k-2$, so the result follows.
\end{proof}

\section{Appendix}

Here we give a different proof of Proposition~\ref{prop:OrderedComplex} using the argument from Proposition~2.14 in \cite{RWW}. The main step is a version of Theorem~2.4 of \cite{GRW} which we give as a lemma below that is of interest in its own right. This gives conditions under which a simplicial map $f\colon Y\to X$ can be homotoped to be injective on individual simplices of some subdivision of $Y$.  We in fact need a relative version of this in which $f$ and the triangulation are kept fixed on a subcomplex $Z\subset Y$.  In this case the best one could hope for is that $f$ is {\em simplexwise injective relative to $Z$}, meaning that either of the following two equivalent conditions is satisfied:

\begin{itemize}
\item If an edge $[v,w]$ has $f(v)=f(w)$, then $[v,w]\in Z$.
\item  For each vertex $v$ of $Y-Z$ we have $f(\lk(v))\subset \lk(f(v))$.
\end{itemize}

\begin{lemma}  Let $X$ be a wCM complex of level $n$, $Y$ a finite simplicial complex of dimension at most $n$,  $f_0\colon Y\to X$ a simplicial map and $Z$ a subcomplex of $Y$.  Then $f_0$ is homotopic fixing $Z$ to a new map $f_1$ that is simplicial with respect to a new triangulation of  \,$Y$ subdividing the old one and
unchanged on $Z$, such that $f_1$ is simplexwise injective relative to $Z$ in the new triangulation of \,$Y$.
\end{lemma}

\begin{proof} 
The proof is by induction on $n$ using a link argument, where the induction starts with the trivial case $n=0$.  Define a simplex $\sigma$ of $Y$ to be {\em bad\/} if for each vertex $v$ of $\sigma$ there is another vertex $w$ of $\sigma$ with $f_0(v)=f_0(w)$.  Our goal is to eliminate all bad simplices that are not contained in $Z$, and in particular all bad edges not contained in $Z$.

If there are any bad simplices not contained in $Z$, let $\sigma$ be one of maximal dimension $k$ (note that $k>0$ since vertices cannot be bad). Since $f_0$ is simplicial we have $f_0(\lk(\sigma))\subset \st(f_0(\sigma))$, but by maximality of $\sigma$ we actually have $f_0(\lk(\sigma))\subset \lk(f_0(\sigma))$. (If $v\in \lk(\sigma)$ maps to $f_0(v)\in f_0(\sigma)$, then $\sigma*v$ is bad, contradicting maximality of $\sigma$.)

Since $\sigma$ is bad, $f_0(\sigma)$ is a simplex of dimension at most $k-1$, hence $\lk(f_0(\sigma))$ is wCM of level $n-k$ since we assumed $X$ is wCM of level $n$.  We assumed also that $Y$ has dimension at most $n$ so $\lk(\sigma)$ has dimension at most $n-k-1$.  Therefore the map $f_0\colon\lk(\sigma)\to \lk(f_0(\sigma))$ is nullhomotopic and we can extend it to a map $g_0\colon b*\lk(\sigma)\to \lk(f_0(\sigma)),$ where $b$ is the barycenter of $\sigma$.  Since $k>0$ we have $n-k<n$ and we can apply induction to deform $g_0$ to a map $g_1$ which agrees with $g_0$ on $\lk\sigma$ and is simplex-wise injective relative to $\lk(\sigma)$ in some subdividion of $b*\lk(\sigma)$ that is unchanged on $\lk(\sigma)$. 
  This homotopy extends over $\st(\sigma)=\sigma*\lk(\sigma)=\bdry\sigma*b*\lk(\sigma)$ by taking the join with the constant homotopy of $f_0\colon \bdry\sigma\to f_0(\sigma)$.  The resulting map $f_1\colon \st(\sigma)\to f_0(\sigma)*\lk(f_0(\sigma))=\st(f_0(\sigma))$ is simplexwise injective relative to $\bdry\sigma *\lk(\sigma)$.  We now have two maps $f_0$ and $f_1$ from $\st(\sigma)$ to $\st(f_0(\sigma))$  which agree on $\bdry\sigma*\lk(\sigma)$.  Since $\st(f_0(\sigma))$ is contractible, these two maps are homotopic by a homotopy which fixes $\bdry\sigma*\lk(\sigma)$.  Extend this homotopy by the constant homotopy outside $\st(\sigma)$.  The resulting map $f_1\colon Y\to X$ now has no bad simplices in the (subdivided) $\st(\sigma)$ except those in $\bdry \sigma *\lk(\sigma)$ that were present before the modification.  The process can now be repeated for other bad simplices of dimension $k$ not contained in $Z$ until they are all eliminated. Then we proceed to $(k-1)$-simplices, etc.  
\end{proof}

\smallskip\noindent {\bf Remark.}  The proof used the connectivity assumptions on links of simplices in $X$ but not that $X$ itself is $(n-1)$-connected. Thus we really only need local conditions on $X$, as one might expect.

\medskip
We now give the alternative proof of Proposition~\ref{prop:OrderedComplex}. Recall that this states that if a simplicial complex $X$ is wCM of level $n$ then its ordered version $\ord{X}$ is $(n-1)$-connected. (It follows that $\ord{X}$ is wCM of level $n$ as well, using the natural extension of this notion to semi-simplicial complexes.)

By induction on $n$ it suffices to show that a map $\bdry D^n\to \ord{X}$ can be extended over $ D^n$.  We first show that any map  $\bdry D^n\to \ord{X}$ is homotopic to a simplicial map.  We cannot directly appeal to the simplicial approximation theorem here because $\ord{X}$ is only semi-simplicial.  (The simplicial approximation theorem does generalize to the semi-simplicial setting as shown in Theorem~5.1 of \cite{RS}, but we do not need the full strength of this.) 

Given $\ord{f}\colon\bdry D^n\to \ord{X},$  let $f\colon \bdry D^n\to X$ be the  composition of $\ord{f}$ with the projection  $p\colon \ord{X}\to X$. Since $X$ is $(n-1)$-connected we can extend  $f$ to $F\colon D^n\to X$. Since $X$ {\em is} a simplicial complex we can use  the simplicial approximation theorem to get a homotopy from $F$ to a map $G$ that is simplicial in a PL triangulation of $D^n$ subdividing any given triangulation.  We assume that this homotopy is constructed as in the standard proof of simplicial approximation, in which case   
 the restriction of the homotopy to $\bdry D^n$ lifts to a homotopy of $\ord{f}$.  This is because the homotopy from $F$ to $G$ has the property that if $F(x)$ lies in the interior of a simplex $\sigma$ of $X$, then $G(x)$ also lies in $\sigma$ (possibly in $\bdry \sigma$) and the homotopy just moves $F(x)$ along the linear path to $G(x)$ in $\sigma$. 
If $x$ lies in $\bdry D^n$ then $\sigma$ lifts to an ordered simplex containing $\ord{f}(x)$  
and the path from $F(x)$ to $G(x)$ also lifts to this ordered simplex.  
 
Since we are free to deform the original map $\ord{f}\colon \bdry D^n\to\ord{X}$ by any homotopy before extending it over $D^n$, we may therefore assume that $\ord{f}$ is simplicial from the start and that we have a simplicial extension $F\colon D^n\to X$ of $f=p\ord{f}$.

We now apply the preceding lemma with $(Y,Z)=(D^n,\bdry D^n)$ to obtain a new map $F\colon D^n\to X$ that is simplexwise injective on a subdivided triangulation of $D^n$, relative to $\bdry D^n$.  This $F$ can be lifted to $\ord{X}$ in the following way.  Choose a total ordering for  the  interior vertices of $D^n.$   
This total ordering gives an ordering on the vertices of each simplex $\tau$ which has no vertices in $\bdry D^n$, and these orderings are compatible with passing to faces.  Since $F$ is injective on each such simplex $\tau$  the ordering on $\tau$ carries over to an ordering of $F(\tau)$ compatible with passing to faces.   This gives a continuous lift $\ord{F}$ of $F$ on interior simplices of $D^n$.  

We already have a lift $\ord{f}$ of $F$ on $\bdry D^n$.  It  remains to lift $F$ on the simplices of $D^n$ that meet $\bdry D^n$ but are not contained in it.  If every such simplex is the join of a boundary simplex $\sigma$ with an interior simplex $\tau$ then the ordering on $f(\sigma)$ given by $\ord{f}$ extends to an ordering on $F(\sigma*\tau)$ by orienting each edge with one end in $f(\sigma)$ and the other in $F(\tau)$ towards $F(\tau)$, i.e., ordering all vertices of $f(\sigma)$ before all vertices of $F(\tau)$.  This works since the lemma guarantees that no such edge is collapsed by $F$ to a single vertex.  

It is possible that  simplices meeting $\bdry D^n$ are  {\em not} joins of boundary and interior simplices, for example an  edge passing through the interior of $D^n$ might have both vertices in  $\bdry D^n$.  
 To avoid this situation,  
note first that the join property is preserved under subdivision. If we start at the very beginning with a triangulation of $D^n$ that has this property, for example by coning off a triangulation of $\bdry D^n$, then the initial simplicial approximation step in the proof gives a subdivision of this triangulation, and applying the lemma produces a further subdivision.

%
%
%
%


\begin{thebibliography}{1}

\bibitem{A} V. I. Arnold,
\newblock Certain topological invariants of algebraic functions, (Russian) 
\newblock {\em Trudy Moskov. Mat. Obsc.} 21 (1970) 27--46. 
\newblock English transl. in {\em Trans. Moscow Math. Soc.} 21 (1970), 30--52. 

\bibitem{A2} V. I. Arnold,
\newblock The cohomology ring of the colored braid group, (Russian) 
\newblock {\em Trudy Moskov. Mat. Obsc.} 5 (1969), 227Ð-231. 
\newblock English transl. in {\em Mathematical Notes} 5 (1969), 138--140.

\bibitem{BBM} Joan Birman,  Nathan Broaddus, and William Menasco, 
\newblock Finite rigid sets and homologically nontrivial spheres in the curve complex of a surface,
\newblock  {\em J. Topol. Anal.} 7 (2015), no. 1, 47--71. 

\bibitem{B} S{\o}ren K. Boldsen,
\newblock Improved homological stability for the mapping class group with integral or twisted coefficients,
\newblock {\em Math. Z.} 270 (2012), 297--329.

\bibitem{Broaddus} Nathan Broaddus,
\newblock Homology of the curve complex and the Steinberg module of the mapping class group,
\newblock {\em Duke Math. J.} 161 (2012), no. 10, 1943--1969. 

\bibitem{Br}
Kenneth S. Brown,
\newblock {\em Cohomology of Groups,}
\newblock Graduate Texts in Mathematics No. 87  (1982) Springer-Verlag, New York.

\bibitem{Charney} R. Charney, 
\newblock Homology stability for $GL_n$ of a Dedekind domain, 
\newblock {\em Invent. Math.} 56 (1980), 1--17.

\bibitem{CD} Ruth Charney and Michael W. Davis,
\newblock Finite $K(\pi,1)$'s for Artin groups,
\newblock  {\em Prospects in Topology,}  Ann. of Math. Stud. 138 (1995), 110--124.

\bibitem{FM}
Benson Farb and Dan Margalit,
\newblock {\em A Primer on Mapping Class Groups,}
\newblock Princeton University Press, 2012.

\bibitem{GRW}
S{\o}ren Galatius and Oscar Randal-Williams,
\newblock Homological stability for moduli spaces of high dimensional manifolds I,
\newblock arXiv:1403.2334

\bibitem{Har1} John L. Harer,
\newblock Stability of the homology of the mapping class groups of orientable surfaces,
\newblock {\em Ann. of Math.} 121 (1985), 215--249.

\bibitem{Har2} John L. Harer,
\newblock The virtual cohomological dimension of the mapping class group of an orientable surface, 
\newblock {\em Invent. Math.} 84(1) (1986), 157--176.

\bibitem{H1}
Allen Hatcher,
\newblock On triangulations of surfaces,
\newblock {\em Topology Appl.}  40 (1991), 189--194. 
\newblock (Updated version at http://www.math.cornell.edu/$\tilde{\ }$hatcher/Papers/TriangSurf.pdf)

\bibitem{H2}
Allen Hatcher, 
\newblock Homological stability for automorphism groups of free groups,
\newblock {\em Comment. Math. Helv.} 70 (1995), 39--62.

\bibitem{HV1}
Allen Hatcher and Karen Vogtmann,
\newblock Cerf theory for graphs,
\newblock  {\em J. London Math. Soc.} 58 (1998), 633--655.

\bibitem{HV2}
Allen Hatcher and Karen Vogtmann,
\newblock  Homology stability for outer automorphism groups of free groups,
\newblock {\em Alg. \& Geom. Topology} 4 (2004), 1253--1272. 

\bibitem{HVW}
Allen Hatcher, Karen Vogtmann, and Nathalie Wahl,
\newblock Erratum to: Homology stability for outer automorphism groups of free groups, 
\newblock {\em Alg. \& Geom. Topology} 6 (2006), 573--579.

\bibitem{HW}
Allen Hatcher and Nathalie Wahl,
\newblock Stabilization for mapping class groups of 3-manifolds,
\newblock {\em Duke Math. J.} 155 (2010), 205--269.

\bibitem{PS}
Andrew Putman and Steven Sam,
\newblock Representation stability and finite linear groups,
\newblock arXiv:1408.3694

\bibitem{Q}
Daniel Quillen,
\newblock Higher Algebraic K-theory I,
\newblock {\em Lecture Notes in Mathematics} 341 (1972), 85--147.

\bibitem{R-W}
Oscar Randall-Williams,
\newblock Resolutions of moduli spaces and homological stability, {\em J. Eur. Math. Soc.} 18 (2016), no. 1, 1--81.

\bibitem{RWW}
Oscar Randall-Williams and Nathalie Wahl,
\newblock Homological stability for automorphism groups,
\newblock arXiv:1409.3541

\bibitem{RS} C. P. Rourke and B. J. Sanderson, 
\newblock $\Delta$-sets I: homotopy theory, 
\newblock {\em Quart. J. Math. Oxford}  22 (1971), 321--338.

\bibitem{vanderKallen} W. van der Kallen, 
\newblock Homology stability for linear groups,
\newblock {\em Invent. Math.} 60 (1980), no. 3, 269--295. 

\bibitem{Vogtmann}K. Vogtmann, 
\newblock Homology stability for $O_{n,n}$, 
\newblock {\em Comm. Algebra}  7 (1979), 9--38.
 
\bibitem{Wagoner}J. B. Wagoner, 
\newblock Stability for homology of the general linear group of a local ring, 
\newblock {\em Topology} 15 (1976), 417--423.

\bibitem{W1}
Nathalie Wahl,
\newblock Homological stability for the mapping class groups of non-orientable surfaces,
\newblock {\em  Invent. Math.} 171 (2008), 389--424.

\bibitem{W2}
Nathalie Wahl,
\newblock Homological stability for mapping class groups of surfaces,
\newblock  {\em Handbook of Moduli, Vol. III,} 
\newblock Advanced Lectures in Mathematics 26 (2012), 547--583. 

\bibitem{Weiss}
Michael Weiss,
\newblock What does the classifying space of a category classify?,
\newblock {\em  Homol. Homot. \& Appl.} 7 (2005), 185--195.



\end{thebibliography}
\end{document}